\documentclass[12pt,twoside,a4paper]{article}

\usepackage{fullpage}
\usepackage{enumerate}
\usepackage{subfigure}
\usepackage{xcolor}

\usepackage[utf8]{inputenc}
\usepackage[T1]{fontenc}
\usepackage[english]{babel}
\usepackage{enumitem}
\usepackage{amsmath,amsfonts,amssymb,graphicx,bbm}
\usepackage{stackrel}

\usepackage{amsthm}
\usepackage{constants}
\newconstantfamily{small}{symbol=c}

\theoremstyle{plain}
\newtheorem{theorem}{Theorem}

\newtheorem{proposition}{Proposition}[section]
\newtheorem{lemma}[proposition]{Lemma}
\newtheorem{corollary}[proposition]{Corollary}

\theoremstyle{definition}
\newtheorem*{definition}{Definition}
\newtheorem{remark}[proposition]{Remark}

\usepackage{mathtools}
\DeclarePairedDelimiter\abs{\lvert}{\rvert}

\newcommand{\Cc}{\mathcal{C}}
\newcommand{\calD}{\mathcal{D}}

\newcommand{\calH}{\mathcal{H}}

\newcommand{\calP}{\mathcal{P}}

\newcommand{\bbE}{\mathbb{E}}

\newcommand{\bbL}{\mathbb{L}}

\newcommand{\bbN}{\mathbb{N}}

\newcommand{\prob}{\mathbb{P}}

\newcommand{\bbR}{\mathbb{R}}

\newcommand{\bbT}{\mathbb{T}}

\newcommand{\bbZ}{\mathbb{Z}}

\newcommand{\Z}{\mathbb{Z}}
\newcommand{\T}{\mathbb T}

\newcommand{\Prob}[1]{\mathbb{P} \left[ #1 \right] }

\newcommand{\R}{\mathbb{R}}

\newcommand{\eps}{\varepsilon}

\def\1{\mathbbm{1}}

\newcommand{\xlra}{\xleftrightarrow}
\newcommand\concel[2]{\ooalign{$\hfil#1\mkern0mu/\hfil$\crcr$#1#2$}}  
\newcommand\nxlra[1]{\mathrel{\mathpalette\concel{\xlra{#1}}}}

\renewcommand{\int}{\mathrm{in}}

\newcommand{\sfe}{\mathsf{e}}
\newcommand{\atrc}{\mathrm{ATRC}}
\newcommand{\betasd}{\beta_{\mathrm{sd}}}

\newcommand{\betac}{\beta_{\mathrm{\scriptscriptstyle c}}}
\newcommand{\norm}[1]{\|#1\|}
\newcommand{\Cov}{\mathrm{Cov}}

\newcommand{\HF}{\mathsf{HF}}
\newcommand{\Spin}{\mathsf{Spin}}

\newcommand{\FK}{\mathsf{FK}}

\newcommand{\ind}[1]{\mathbbm{1}_{\{#1\}}}

\usepackage{verbatim}

\usepackage{tabularx}

\usepackage{hyperref}

\title{Phase diagram of the Ashkin--Teller model}
\date{\today}

\author{
	Yacine Aoun
		\thanks{
		Section de Mathématiques, 
		Université de Genève, 
		CH-1211 Genève, 
		Switzerland
		\url{Yacine.Aoun@unige.ch}
		}
	\and
	Moritz Dober
		\thanks{
		Faculty of Mathematics, 
		University of Vienna, 
		Oskar-Morgenstern-Platz 1, 
		A-1090 Vienna, Austria.
		\url{moritz.dober@univie.ac.at}
		}
	\and
	Alexander Glazman
    	\thanks{Department of Mathematics,
    	University of Innsbruck,
	    Technikerstr. 13,
	    A-6020 Innsbruck, 
	    Austria.
	    \url{alexander.glazman@uibk.ac.at}}	    
	}

\date{\today}

\begin{document}

\maketitle

\begin{abstract}
The Ashkin--Teller model is a pair of interacting Ising models and has two parameters: $J$ is a coupling constant in the Ising models and~$U$ describes the strength of the interaction between them.
In the ferromagnetic case~$J,U>0$ on the square lattice, we establish a complete phase diagram conjectured in physics in 1970s (by Kadanoff and Wegner, Wu and Lin, Baxter and others): when~$J<U$, the transitions for the Ising spins and their products occur at two distinct curves that are dual to each other; when~$J\geq U$, both transitions occur at the self-dual curve.
All transitions are shown to be sharp using the OSSS inequality.

We use a finite-criterion argument and continuity to extend the result of Peled and the third author~\cite{GlaPel19} from a self-dual point to its neighborhood.
Our proofs go through the random-cluster representation of the Ashkin--Teller model introduced by Chayes--Machta and Pfister--Velenik and we rely on couplings to FK-percolation.
\end{abstract}

\section{Introduction}
\label{sec:intro}

The Ashkin--Teller (AT) model is named after two physicists who introduced it in 1943~\cite{AshTel43} and can be viewed as a pair of interacting Ising models.
For a finite subgraph~$\Omega=(V,E)$ of~$\bbZ^2$, the AT model is supported on pairs of spin configurations~$(\tau,\tau') \in \{\pm1\}^V \times \{\pm1\}^V$ and the distribution is defined by
\begin{equation}\label{eq:def-at}
	{\sf AT}_{\Omega,J_{\tau},J_{\tau'},U}(\tau,\tau') = \frac{1}{Z}\cdot \exp \left[\sum_{uv\in E} J_{\tau}\tau_u\tau_v + J_{\tau'}\tau_u'\tau_v' + U\tau_u\tau_u'\tau_v\tau_v' \right],
\end{equation}
where~$J_{\tau},J_{\tau'},U$ are real parameters and~$Z = Z(\Omega,J_{\tau},J_{\tau'},U)$ is the unique constant (called {\em partition function}) that renders the above a probability measure. 

In the current article, we consider the ferromagnetic symmetric (or isotropic) case
\[
	J=J_{\tau}=J_{\tau'}\geq 0, \text{ and } U\geq 0
\]
and denote the measure by~${\sf AT}_{\Omega,J,U}$.

Important particular cases: $U=0$ gives two independent Ising models;for $J=0$, $\tau$ reduces to a Bernoulli site percolation with parameter $1/2$, and $\tau\tau'$ to an Ising model, independent of each other; the line~$U=J$ corresponds to the 4-Potts model.
These models are very well-studied and their phase diagram is known; see~\cite{FriVel17,Dum17a} for excellent surveys. 
Henceforth in this article we assume that $J,U>0$.
A key observation in the analysis of the general AT model is its relation to the six-vertex model~\cite{Fan72}.
This gives a {\em non-staggered} six-vertex model (i.e. with shift invariant local weights) only at the {\em self-dual line} of the AT model: it was found in~\cite{MitSte71} and is described by the equation
\[
	\label{eq:self-dual}\tag{SD}
	\sinh (2J) = e^{-2U}.
\]
Outside of this line, the corresponding six-vertex model is staggered and thus the seminal Baxter's solution~\cite{Bax71} does not apply.
Kadanoff and Wegner~\cite{KadWeg71, Weg72}, Wu and Lin~\cite{WuLin74}, and others conjectured that, when~$J<U$, there are two {\em distinct} transition lines in the AT model: one for correlations of spins~$\tau$ (or~$\tau'$) and the other for correlations of products~$\tau\tau'$.
In the current article, we prove this conjecture and establish a complete phase diagram of the AT model in the ferromagnetic regime. 

It will be convenient to state the results in infinite volume and to consider also plus boundary conditions.
Denote by~$\partial \Omega$ the set of boundary vertices of~$\Omega$~-- these are all vertices in~$\Omega$ that are adjacent to at least one vertex in $\bbZ^2\setminus\Omega$.
We define the measure with plus boundary conditions by conditioning all boundary vertices to have spin plus in~$\tau$ and in~$\tau'$:
\[
	{\sf AT}_{\Omega,J,U}^{+,+} := {\sf AT}_{\Omega,J,U}(\cdot \, | \, \tau_{|\partial \Omega} \equiv \tau_{|\partial \Omega}' \equiv 1).
\] 
Expectations with respect to the AT measures are denoted by brackets: 
\[
	\langle \cdot \rangle_{\Omega,J,U}:= \bbE_{\Omega,J, U}[\cdot ] \hspace{10mm} \text{and} \hspace{10mm}
	\langle \cdot \rangle_{\Omega,J,U}^{+,+}:= \bbE_{\Omega,J, U}^{+,+}[\cdot ].
\]

The correlations satisfy the Griffiths--Kelly--Sherman (GKS) inequality~\cite{KelShe68}, which states that for any $A,B,C,D\subset V$, one has 
\begin{equation}\label{eq:gks}\tag{GKS}
\langle \tau_{A}\cdot\tau'_{B}\cdot\tau_{C}\cdot\tau'_{D}\rangle _{\Omega,J,U}\geq 
\langle \tau_{A}\cdot\tau'_{B}\rangle _{\Omega,J,U}\langle\tau_{C}\cdot\tau'_{D}\rangle _{\Omega,J,U},
\end{equation}
where~$\tau_A:=\prod_{u\in A}\tau_u$ and~$\tau_B':=\prod_{v\in B}\tau_v'$. This in particular implies that any~$A,B\subset V$,
\[
	\langle \tau_A \cdot \tau_B' \rangle_{\Omega,\beta J,\beta U} \text{ and }
	\langle \tau_A \cdot \tau_B' \rangle_{\Omega,\beta J,\beta U}^{+,+} \text{ are increasing in } \beta>0,
\]

A standard application of the GKS inequality implies that the weak limits over~$\Omega_n\nearrow \bbZ^2$ exist and do not depend on~$\{\Omega_n\}$:
\[
	{\sf AT}_{J,U} := \lim_{n\to \infty} {\sf AT}_{\Omega_n,J,U} \hspace{10mm} \text{and} \hspace{10mm}
	{\sf AT}_{J,U}^{+,+} := \lim_{n\to \infty} {\sf AT}_{\Omega_n,J,U}^{+,+}.
\]
Similarly to finite-volume measures, we denote by~$\langle \cdot \rangle_{J,U}$ and~$\langle \cdot \rangle_{J,U}^{+,+}$ the correlation functions with respect to~${\sf AT}_{J,U}$ and~${\sf AT}_{J,U}^{+,+}$.
It is standard (e.g. can be shown by comparing to the Ising model) that the AT model undergoes a phase transition in terms of correlations of~$\tau$ and those of~$\tau\tau'$.
Moreover, a general OSSS inequality~\cite{DumRaoTas19} can be used to show that both transitions are sharp (see Appendix). 
That is, for each pair~$J,U$, there exist~$\beta_c^{\tau}, \beta_c^{\tau\tau'}, c, C > 0$, such that 
\[
	\langle\tau_0\tau_x\rangle_{\beta J, \beta U}^{+,+}
	\begin{cases} 
		\leq e^{-c\cdot \parallel x\parallel} & \text{if } \beta<\beta_c^{\tau} \\
		\geq C &  \text{if } \beta>\beta_c^{\tau} \\
	\end{cases}
	,
	\qquad
	\langle \tau_0\tau_0'\tau_x\tau_x'\rangle_{\beta J, \beta U}^{+,+}
	\begin{cases} 
		\leq e^{-c\cdot \parallel x\parallel} & \text{if } \beta<\beta_c^{\tau\tau'} \\
		\geq C &  \text{if } \beta>\beta_c^{\tau\tau'} \\
	\end{cases} \, .
\]
Symmetry between $\tau$ and $\tau'$ and the correlation inequalities~\eqref{eq:gks} imply directly that
\begin{equation}\label{eq:easy-betac-ineq}
	\betac^{\tau}\geq\betac^{\tau\tau'}.
\end{equation}
We define the transition points with respect to the free measure in the similar way, and we denote them by $\betac^{\tau,f}$ and $\betac^{\tau\tau',f}$.

There exists a unique~$\beta$, for which~$(\beta J, \beta U)$ is on the line~\eqref{eq:self-dual}.
Denote it by~$\betasd=\betasd(J,U)$.
The following theorem states our main result:

\begin{theorem}\label{thm:splitting}
	Let~$0<J<U$. Then, $\beta_c^{\tau} >\betasd > \beta_c^{\tau\tau'}$.
\end{theorem}
This was previously shown when~$2J<U$ using a direct comparison to the Ising model~\cite{Pfi82}.
In addition, in the perturbative regime when~$J/U$ is big enough, the critical exponents associated to both phase transitions have been shown to be the same as for the Ising model~\cite{GiuMas05}. It is expected that the critical exponents vary continuously in the whole regime $J<U$, and that the critical exponents are the same as for the Ising model when $J\geq U$. We refer to~\cite{DelGri04} for a survey on the physics literature on the critical behaviour of the AT model, as well as predictions on critical exponents using the quantum filed theory. Recently, Peled and the third author have proven that spins~$\tau$ (or~$\tau'$) and the products~$\tau\tau'$ exhibit qualitatively different behavior at the self-dual line when~$U>J$~\cite{GlaPel19}:
products~$\tau\tau'$ are ordered, while~$\tau$ (and $\tau'$) exhibits exponential decay of correlations.
We derive Theorem~\ref{thm:splitting} by extending this statement to an open neighborhood of the self-dual line when~$J<U$.
The continuity ideas do not apply directly, since the rate of decay of correlations might, a priori, depend on infinitely many spins.
To circumvent this problem, we establish exponential decay in finite volume:

\begin{proposition}\label{thm:finite_expo_decay_at_sd}
	Fix $0<J<U$ that satisfy~$\sinh 2J = e^{-2U}$. 
	Then, there exists $c:=c(J,U)>0$ such that
	\begin{equation*}
		\langle\tau_{0}\rangle^{+,+}_{[-n,n]^2, J,U}\leq e^{-cn}.
	\end{equation*}
\end{proposition}

\begin{figure}
	\begin{center}		
		\includegraphics[scale=0.8]{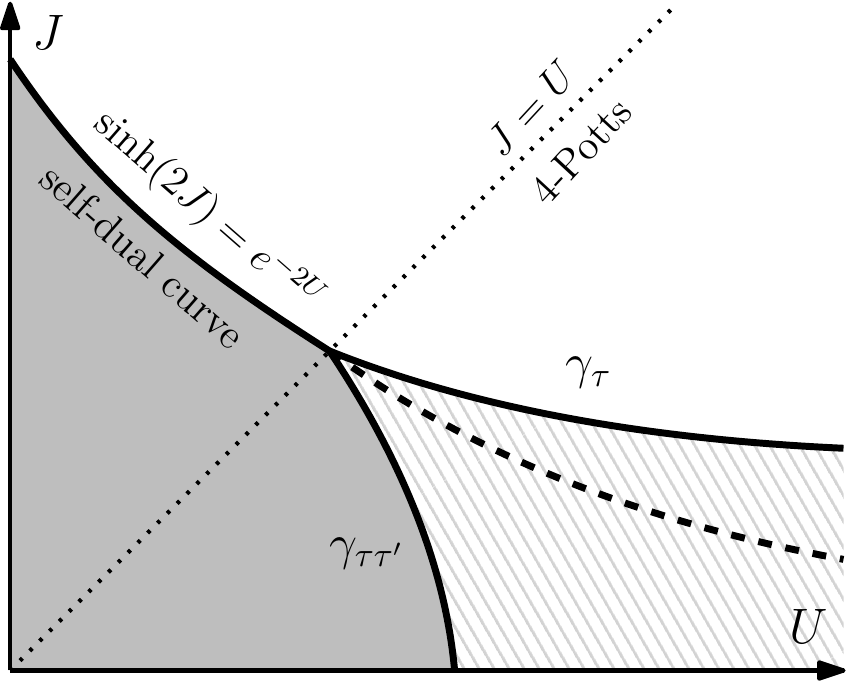}\hspace{20mm}
		\includegraphics[scale=0.9,page=3]{BKW.pdf}
	\end{center}		
	\caption{{\em Left:} Phase diagram of the Ashkin--Teller model: when~$J\geq U$, transitions for~$\tau$ and~$\tau\tau'$ occur at the self-dual curve (Theorem~\ref{thm:sd-is-critical}) and when~$J < U$, the transition occurs at two distinct curves~$\gamma_\tau$ and~$\gamma_{\tau\tau'}$ dual to each other (Theorems~\ref{thm:splitting} and~\ref{thm:duality}).
	There are three regimes: disorder in~$\tau$ and in~$\tau\tau'$ (gray), order in~$\tau$ and in~$\tau'$ (white), disorder in~$\tau$ and order in~$\tau\tau'$ (dashed gray).
	{\em Right:} Domain $\Omega$ (in bold black) on $\bbL$ and its dual $\Omega^*$ (in gray). Notice that~$\Omega^*$ is not a domain on $\bbL^{*}$. The even domain $\calD_{\Omega}$ (dashed) on~$\bbZ^2$.}
	\label{fig:at-diagram}
\end{figure}

Compared to~\cite{GlaPel19}, the exponential decay is proven in the finite-volume and under the largest boundary conditions.
This is crucial for applying the so-called ``finite criterion'' (or $\varphi_\beta(S)$) argument~\cite{Sim80,Lie80,DumTas15,DumTas15a}, since the Simon--Lieb inequality is not available.
This argument, as well as the proof of Proposition~\ref{thm:finite_expo_decay_at_sd}, use the random-cluster representation of the AT model (that we call ATRC) introduced by Chayes--Machta~\cite{ChaMac97} and Pfister--Velenik~\cite{PfiVel97}.
As in the seminal Edwards--Sokal coupling for the Potts model~\cite{EdwSok88}, connectivities in the ATRC describe correlations in the AT model.

Other ingredients in the proof of Proposition~\ref{thm:finite_expo_decay_at_sd} are couplings between the AT and the six-vertex models~\cite{Fan72,Weg72} and between the latter and FK-percolation~\cite{BaxKelWu76}.
These two couplings were composed for the first time in the work of Peled and the third author~\cite{GlaPel19}.
We also use $\bbT$-circuits introduced in~\cite{GlaPel19} to apply the non-existence theorem~\cite{She05,DumRaoTas19}.

The next result states that the transition lines are dual to each other (see Fig.~\ref{fig:at-diagram}) and that the critical points for the measures under the free and plus boundary conditions coincide. In order to make a precise statement, we define the critical curve 
\begin{equation*}
\gamma_{\tau}:=\lbrace (J,U)\in\bbR^{2} : 0<J<U,\text{ }\betac^{\tau}(J,U)=1\rbrace.
\end{equation*}
We define $\gamma_{\tau\tau'}$ in a similar way. Given a pair of parameters $(J,U)$, we define the dual set of parameters $(J^{*},U^{*})$ as the unique solutions to the following equations
	\begin{equation}\label{eq:duality_1}
		\frac{e^{-2J+2U}-1}{e^{-2J^*+2U^*}-1} = e^{2U}\sinh (2J) = \left[e^{2U^*}\sinh (2J^*)\right]^{-1}.
	\end{equation}
Note that this duality relation is an involution. We refer the reader to Subsection~\ref{sec:rc-at} and Lemma~\ref{lemma:duality_relation} for more details on the duality in the ATRC model. 

\begin{theorem}\label{thm:duality}
Fix~$0<J<U$. Then, the following holds:
\begin{enumerate}[label=(\roman*)]
	\item $\betac^{\tau}=\betac^{\tau,f}$ and $\betac^{\tau\tau'}=\betac^{\tau\tau',f}$;
	\item $\gamma_{\tau}$ and $\gamma_{\tau\tau'}$ are dual in the following sense: $(J,U)\in\gamma_{\tau}$ if only if $(J^{*},U^{*})\in\gamma_{\tau\tau'}$.
\end{enumerate}
\end{theorem}

The next theorem states that when~$J\geq U$, both transitions in~$\tau$ and in~$\tau\tau'$ occur at the self-dual line.

\begin{theorem}\label{thm:sd-is-critical}
	Let~$J\geq U > 0$. Then, the following holds:
\begin{enumerate}[label=(\roman*)]
	\item $\betac^{\tau}=\betac^{\tau,f}$ and $\betac^{\tau\tau'}=\betac^{\tau\tau',f}$;
	\item $\betac^{\tau} = \betac^{\tau\tau'} =\betasd$.
\end{enumerate}	
\end{theorem}

General approach~\cite{DumRaoTas19} gives sharpness under plus boundary conditions and equality of the transition points~$\betac^{\tau}= \betac^{\tau\tau'}=:\betac$.
By standard duality arguments, one deduces~$\betac \leq \betasd$.
The bound~$\betac \geq \betasd$ follows from Zhang-type arguments provided the transition points for the free and monochromatic measures coincide.
The latter can be shown by applying the classical FK-percolation argument to the marginals of the ATRC.

\paragraph{Organisation of the article.}
Sections~\ref{section:proof_splitting}--\ref{sec:exp-dec-finite} treat the case $J<U$:
in Section~\ref{section:proof_splitting}, we introduce the random-cluster representation of the AT model (ATRC) and derive Theorems~\ref{thm:splitting} and~\ref{thm:duality} from Proposition~\ref{thm:finite_expo_decay_at_sd}; Sections~\ref{sec:couplings}--\ref{sec:exp-dec-finite} are dedicated to proving Proposition~\ref{thm:finite_expo_decay_at_sd}.
In Section~\ref{sec:couplings}, we describe the six-vertex and FK-percolation models and give their background, including their relation to the AT model.
In Section~\ref{sec:exp-rel-0-1}, we show that~$\tau$ exhibits exponential decay of correlations in finite volume under the boundary conditions~$\tau=\tau'$.
In Section~\ref{sec:no-inf-cluster-1-1}, we show that~$\tau$ exhibits no ordering under~${\sf AT}_{J,U}^{+,+}$.
In Section~\ref{sec:exp-dec-finite}, we derive Proposition~\ref{thm:finite_expo_decay_at_sd}. Section~\ref{sec:0-le-u-le-j} deals with the case $J\geq U$: we introduce the ATRC model and prove Theorem~\ref{thm:sd-is-critical}.
Appendices provide details regarding sharpness for the AT (\ref{sec:sharpness}), exponential relaxation for FK-percolation (\ref{sec:exp-rel-fk}), stochastic ordering of the ATRC with respect to its local weights (\ref{sec:proof-stoch-dom-atrc}) and uniqueness of the infinite-volume ATRC measure (\ref{sec:dense_subsets}).

\paragraph{Acknowledgements}
We would like to thank Aran Raoufi for sharing his notes on sharpness from~2017 and Yvan Velenik for directing us towards~\cite[Appendix]{CamIofVel08}.
We are also grateful to Ioan Manolescu and Sebastien Ott for many fruitful discussions.
Part of this work was done during the visits of YA, MD and AG: we would like to thank the Universities of Geneva, Innsbruck and Vienna for their hospitality.

The work of AG and MD is supported by the Austrian Science Fund grant P34713. YA is supported by the Swiss NSF grant 200021\_200422 and is member of NCCR SwissMAP.

\section{From Proposition~\ref{thm:finite_expo_decay_at_sd} to Theorems~\ref{thm:splitting} and~\ref{thm:duality}}\label{section:proof_splitting}
\label{sec:proof-thm1}

From now on, we will consider the AT model on a rotated square lattice that we denote by $\bbL$: its vertex set is~$\{ (x,y)\in\bbZ^{2} \colon x+y \text{ is even}\}$ and edges connect~$(x,y)$ to~$(x\pm 1, y\pm 1)$, see Figure~\ref{fig:at-diagram}. This is more convenient for the coupling with the six-vertex model (Section~\ref{sec:couplings}).

In this section, we fix $J < U$ and drop them from the notation. In particular, we write~${\sf AT}_{\Omega,\beta}$ for the measure~${\sf AT}_{\Omega,\beta J, \beta U}$.

We start by defining the random-cluster representation of the AT model (ATRC) introduced by Chayes--Machta~\cite{ChaMac97} and Pfister--Velenik~\cite{PfiVel97}.
Using a $\varphi_{\beta}(S)$ argument, we prove that~$(\beta_c^{\tau\tau'}J,\beta_c^{\tau\tau'}U)$ is strictly above the self-dual line.
By duality, this implies that~$(\beta_c^{\tau}J,\beta_c^{\tau}U)$ is strictly below the self-dual line which concludes the proof.

\subsection{ATRC: defintion and basic properties}
\label{sec:rc-at}

The ATRC is reminiscent of the Edwards--Sokal~\cite{EdwSok88} coupling between FK-percolation and the Potts model.
Since the AT model is supported on a pair of spin configurations, the ATRC is supported on a pair of bond percolation configurations.

\paragraph{Percolation configurations.}
For a finite subgraph~$\Omega \subset \bbL$, the sets of its vertices and edges are denoted by~$V_\Omega$ and~$E_\Omega$, respectively. 
We view~$\omega \in \{0,1\}^{E_\Omega}$ as a percolation configuration: we say that~$e$ is {\em open} in $\omega$ if~$\omega(e) = 1$, and otherwise~$e$ is {\em closed}. 
We identify~$\omega$ with a spanning subgraph of $\Omega$ and edges that are open in~$\omega$.
Define~$\vert\omega\vert$ as the number of edges in $\omega$. 
Boundary conditions for~$\omega$ are given by a partition $\eta$ of~$\partial \Omega$.
We define~$k^{\eta}(\omega)$ as the number of connected components in~$\omega$ when all vertices belonging to the same element of partition in $\eta$ are identified.
Two important special cases: $1$ denotes {\em wired} b.c. given by a trivial partition consisting of one element~$\partial\Omega$; $0$ denotes {\em free} b.c. given by a partition of~$\partial\Omega$ into singletons.

\paragraph{Definition of ATRC.}
A configuration of the ATRC model on~$\Omega$ is a pair~$(\omega_{\tau},\omega_{\tau\tau'})$ of percolation configurations on edges of~$\Omega$.
Formally, the ATRC measure is supported on
$(\omega_{\tau},\omega_{\tau\tau'})\in \{0,1\}^{E_\Omega}\times \{0,1\}^{E_\Omega}$.
For~$\beta > 0$ and partitions~$\eta_\tau,\eta_{\tau\tau'}$ of~$\partial\Omega$, the ATRC measure is defined by
\begin{equation}\label{eq:definition_ATRC}
	\atrc^{\eta_\tau,\eta_{\tau\tau'}}_{\Omega,\beta}(\omega_\tau, \omega_{\tau\tau'})
	=\tfrac{1}{Z}\cdot 
	2^{k^{\eta_{\tau}}(\omega_{\tau})+k^{\eta_{\tau\tau'}}(\omega_{\tau\tau'})}\prod_{e\in E}a(\omega_\tau(e),\omega_{\tau\tau'}(e)),
\end{equation}
where~$Z=Z(\Omega,\beta,J,U,\eta_\tau,\eta_{\tau\tau'})$ is a normalizing constant and
\begin{equation}\label{eq:atrc-weights}
	a(0,0):= \sfe^{-2\beta (J+U)}, \,\, a(1,0):=0, \,
	a(0,1):=\sfe^{-4\beta J}- \sfe^{-2\beta(J+ U)}, \, a(1,1):= 1-\sfe^{-4\beta J}.
\end{equation}
Since~$J < U$, we have $a(i,j)\geq 0$ for all $i,j\in\lbrace 0,1\rbrace$. We will also use the notation $\atrc_{\Omega,J,U}^{\eta_\tau,\eta_{\tau\tau'}}$ for the measure with $\beta=1$.

It will be useful to express the measure as 
\begin{equation}\label{eq:atrc-repr.}
	\atrc^{\eta_\tau,\eta_{\tau\tau'}}_{\Omega,\beta}(\omega_\tau, \omega_{\tau\tau'})
	\propto 
	{\rm w}_\tau^{\abs{\omega_\tau}}\, {\rm w}_{\tau\tau'}^{\abs{\omega_{\tau\tau'}\setminus\omega_\tau}}\,
	2^{k^{\eta_{\tau}}(\omega_{\tau})+k^{\eta_{\tau\tau'}}(\omega_{\tau\tau'})}\,\ind{\omega_\tau\subseteq\omega_{\tau\tau'}},
\end{equation}
where
\begin{equation}\label{eq:atrc-weights-2}
{\rm w}_\tau=e^{2\beta U}(e^{2\beta J}-e^{-2\beta J})\quad\text{and}\quad {\rm w}_{\tau\tau'}=e^{2\beta(U-J)}-1.
\end{equation}
In this context, we will refer to the measure as $\atrc_{\Omega,{\rm w}_\tau,{\rm w}_{\tau\tau'}}^{\eta_\tau,\eta_{\tau\tau'}}$. In Section~\ref{sec:modified-atrc}, we will encounter a version of this measure with non-homogeneous weights.

\begin{remark}
	The representation can be extended to $J\geq U$~\cite{PfiVel97}, see Section~\ref{sec:atrc-0uj}.
\end{remark}

There are four special types of boundary conditions given by free/wired~$\eta_\tau$ and free/wired~$\eta_{\tau\tau'}$:
$\atrc_{\Omega,\beta}^{1,1}$ (both wired), $\atrc_{\Omega,\beta}^{0,0}$ (both free), $\atrc_{\Omega,\beta}^{1,0}$ (wired for~$\omega_\tau$, free for~$\omega_{\tau\tau'}$), $\atrc_{\Omega,\beta}^{0,1}$ (free for~$\omega_\tau$, wired for~$\omega_{\tau\tau'}$).

\paragraph{Coupling between ATRC and AT.}
For $X,Y\subset\bbL$ and a percolation configuration~$\omega \in \{0,1\}^{E_\Omega}$, we define $X\xlra{\omega} Y$ as an event that~$X$ and $Y$ are linked by a path of open edges in~$\omega$.
If~$X=\{x\}$ and~$Y=\{y\}$, we simply write~$x \xlra{\omega} y$.
We also use the notation $x\xlra{\omega}\infty$ for the event of~$x$ belonging to an infinite connected component of~$\omega$.	

The key property of the ATRC is that connectivities in it describe correlations in the AT model~\cite{PfiVel97}:
for~$\beta> 0$ and any finite subgraph~$\Omega\subset \bbL$ containing~$0$,
\begin{equation}\label{eq:coupling}
	\langle\tau_{0}\rangle^{+,+}_{\Omega,\beta}=\atrc^{1,1}_{\Omega,\beta}(0\xlra{\omega_\tau}\partial\Omega),
	\qquad
	\langle\tau_{0}\tau'_{0}\rangle^{+,+}_{\Omega,\beta}=\atrc^{1,1}_{\Omega,\beta}(0\overset{\omega_{\tau\tau'}}{\longleftrightarrow}\partial\Omega).
\end{equation}
We omit the proof as it is straightforward and similar to the classical Edwards--Sokal coupling; see~\cite[Proposition~3.1]{PfiVel97} for details.

\paragraph{Positive correlations and monotonicity.}
We first introduce the notion of \emph{stochastic domination} and \emph{positive association}. Given a partially ordered set $\mathcal{P}$ and a real-valued function $f$ on $\calP$, $f$ is said to be \emph{increasing} if for any $\omega,\omega'\in\calP$ with $\omega\leq\omega'$, one has $f(\omega)\leq f(\omega')$. A subset $A\subseteq\mathcal{P}$ is then called \emph{increasing} if its indicator $\ind{A}$ is increasing with respect to $\mathcal{P}$. Given two probability measures $\mu$ and $\nu$ on $\mathcal{P}$ equipped with some $\sigma$-algebra $\mathcal{A}$, we say that $\mu$ is stochastically dominated by $\nu$ (or $\nu$ stochastically dominates $\mu$), and write $\mu\leq_{\mathrm{st}}\nu$ (or $\nu\geq_{\mathrm{st}}\mu$), if for every increasing event $A\in\mathcal{A}$, we have $\mu(A)\leq\nu(A)$. Moreover, $\mu$ is said to be positively associated or to satisfy the FKG inequality if, for all increasing non-negative functions $f$ and $g$, we have
\begin{equation}
	\label{eq:fkg-general}
	\mu(f\cdot g)\geq \mu(f)\mu(g).
\end{equation}

We introduce a natural partial order on pairs of percolation configurations: we say that $(\omega_{\tau},\omega_{\tau\tau'}) \geq (\tilde{\omega}_{\tau},\tilde{\omega}_{\tau\tau'})$ if and only if, for every edge~$e$, we have~$\omega_\tau(e)\geq \tilde{\omega}_\tau(e)$ and~$\omega_{\tau\tau'}(e)\geq \tilde{\omega}_{\tau\tau'}(e)$.
By~\cite[Proposition~4.1]{PfiVel97} (and its proof), the measures $\atrc_{\Omega,\beta}^{\eta_\tau,\eta_{\tau\tau'}}$ are positively associated for any~$\beta> 0$ and any boundary conditions~$\eta_\tau,\eta_{\tau\tau'}$: for any increasing events~$A,B$, one has
\[
	\label{eq:fkg} \tag{FKG}
	\atrc_{\Omega,\beta}^{\eta_\tau,\eta_{\tau\tau'}}(A\cap B) \geq 
	\atrc_{\Omega,\beta}^{\eta_\tau,\eta_{\tau\tau'}}(A)\cdot 
	\atrc_{\Omega,\beta}^{\eta_\tau,\eta_{\tau\tau'}}(B).
\]
This can be used to compare different boundary conditions.
For two partitions~$\eta$ and~$\tilde{\eta}$ of~$\partial\Omega$, we say that~$\eta \geq \tilde{\eta}$ if any two vertices belonging to the same element of~$\tilde{\eta}$ also belong to the same element of~$\eta$.
Then, for any~$\beta> 0$, and any boundary conditions such that $\eta_\tau\geq \tilde{\eta}_\tau$ and~$\eta_{\tau\tau'}\geq \tilde{\eta}_{\tau\tau'}$,
\[
	\label{eq:cbc} \tag{CBC}
	\atrc_{\Omega,\beta}^{\eta_\tau,\eta_{\tau\tau'}} \geq_{\rm st} \atrc_{\Omega,\beta}^{\tilde{\eta}_\tau,\tilde{\eta}_{\tau\tau'}}.
\]	
The Holley criterion \cite{Hol74} also allows to show stochastic ordering of the measures in the parameter $\beta$: if $\beta_1\geq\beta_2$, then, for any boundary conditions $\eta_\tau,\eta_{\tau\tau'}$,
\[
	\label{eq:mon} \tag{MON}
	\atrc_{\Omega,\beta_1}^{\eta_\tau,\eta_{\tau\tau'}} \geq_{\rm st} \atrc_{\Omega,\beta_2}^{\eta_\tau,\eta_{\tau\tau'}}.
\]
See~\cite[Lemma 11.14]{Gri06} for a proof. 
In fact, this proof gives a little more. 
Indeed, it consists of checking inequalities for quantities that are continuous functions of $(\beta_i J,\beta_i U)$ and the inequalities are strict when $\beta_1>\beta_2$.
This implies that the ATRC measure with parameters in a small neighbourhood of~$(\beta_1 J, \beta_1 U)$ dominates that with parameters in a small neighbourhood of~$(\beta_2 J, \beta_2 U)$.
More precisely, for~$(x,y)\in \bbR^2$, define $B_{r}(x,y)$ as the Euclidian ball of radius $r$ centred at $(x,y)$. 
If $\beta_1>\beta_2$, then there exists $\eps>0$ such that for any $(J_1,U_1)\in B_\eps(\beta_1 J,\beta_1 U)$ and $(J_2,U_2)\in B_\eps(\beta_2 J,\beta_2 U)$,
\[
	\label{eq:strong-mon} \tag{MON+}
	\atrc_{\Omega,J_1,U_1}^{\eta_\tau,\eta_{\tau\tau'}} \geq_{\rm st} \atrc_{\Omega,J_2,U_2}^{\eta_\tau,\eta_{\tau\tau'}}.
\]
This extension will be useful for our proof of Theorem~\ref{thm:duality}.

\paragraph{Domain Markov property.}
As in the standard FK-percolation, one can interpret a configuration outside of a subdomain as boundary conditions.
Indeed, let~$\Omega \subset \Delta$ be two finite subgraphs of~$\bbL$ and~$\xi\in \{0,1\}^{E_\Delta\setminus E_\Omega}$ a percolation configuration on~$\Delta\setminus \Omega$.
Given boundary conditions~$\eta$ on~$\Delta$, define a partition~$\eta\cup\xi$ of~$\partial\Omega$ by first identifying vertices belonging to the same element of~$\eta$ and then identifying vertices belonging to the same cluster of~$\xi$.
Then, the following domain Markov property holds:
\[
	\label{eq:dmp}\tag{DMP}
	\atrc_{\Delta, \beta}^{\eta_\tau,\eta_{\tau\tau'}}( \cdot \,|\, 
	(\omega_\tau,\omega_{\tau\tau'})_{|\Delta\setminus\Omega} = (\xi_\tau,\xi_{\tau\tau'})_{|\Delta\setminus\Omega})=
	\atrc_{\Omega, \beta}^{\eta_\tau\cup \xi_\tau,\eta_{\tau\tau'}\cup \xi_{\tau\tau'}}( \cdot ).
\]
Thus, by~\eqref{eq:cbc}, for any increasing sequence of subgraphs~$\Omega_k \nearrow \bbL$, the measures~$\atrc_{\Omega_k, \beta}^{1,1}$ form a stochastically decreasing sequence.
Thus, the weak (or local) limit exists and is unique, by standard arguments.
Denote it by~$\atrc_\beta^{1,1}$.
Define~$\atrc_\beta^{0,0}$ analogously. We write~$\atrc_{J,U}^{1,1}$ and~$\atrc_{J,U}^{0,0}$ for the corresponding measures with $\beta=1$. 

\paragraph{Dual ATRC.} 
Define the dual lattice $\bbL^{*}:=\bbL+(1,0)$.
For each edge~$e$ of~$\bbL$, there is a unique edge of~$\bbL^{*}$ that intersects it: call this edge dual to~$e$ and denote it by~$e^*$. Denote by $E_\bbL$ and $E_{\bbL^*}$ the sets of edges of $\bbL$ and $\bbL^*$, respectively. Given a percolation configuration $\omega\in\{0,1\}^{E_\bbL}$, we define its dual configuration $\omega^*\in\{0,1\}^{E_{\bbL^*}}$ by setting
\[
	\omega^*(e^*):=1-\omega(e).
\]
For an ATRC configuration~$(\omega_\tau,\omega_{\tau\tau'})\in \{0,1\}^{E_\bbL}\times \{0,1\}^{E_\bbL}$, we define its dual~$(\hat{\omega}_\tau,\hat{\omega}_{\tau\tau'})\in \{0,1\}^{E_{\bbL^*}}\times \{0,1\}^{E_{\bbL^*}}$ in the following way,
\begin{equation}\label{eq:dualrelation}
	\hat{\omega}_\tau :=\omega_{\tau\tau'}^* \hspace{10mm} \text{and} \hspace{10mm} 
	\hat{\omega}_{\tau\tau'} := \omega_\tau^*.
\end{equation}
We want to emphasize that we are not considering two standard dual percolation configurations but we also swap the order of~$\tau$-edges and~$\tau\tau'$-edges.

The measures $\atrc_{J,U}^{0,0}=:\atrc_{\bbL,J,U}^{0,0}$ and $\atrc_{J,U}^{1,1}=:\atrc_{\bbL,J,U}^{1,1}$ on $\bbL$ can be defined on $\bbL^*$ in the same manner, and we denote them by~$\atrc_{\bbL^*,J,U}^{0,0}$ and~$\atrc_{\bbL^*,J,U}^{1,1}$, respectively. 
Recall the mapping~$(J,U)\mapsto (J^*,U^*)$ defined by~\eqref{eq:duality_1} and note its properties: it is continuous, an involution, identity on the self-dual line~\eqref{eq:self-dual}, sends every point above~\eqref{eq:self-dual} to a point below~\eqref{eq:self-dual}.
The pushforward of the ATRC measure under the duality transformation is also an ATRC measure with the dual parameters:

\begin{lemma}[Prop 3.2 in \cite{PfiVel97}]\label{lemma:duality_relation}
	Let~$0< J < U$.
	Let $(\omega_\tau,\omega_{\tau\tau'})$ be distributed according to $\atrc^{1,1}_{\bbL,J, U}$.
	Then, the distribution of~$(\hat{\omega}_\tau,\hat{\omega}_{\tau\tau'})$ is given by~$\atrc^{0,0}_{\bbL^{*},J^*,U^*}$.
\end{lemma}

\subsection{Proof of Theorem~\ref{thm:duality}}
	\label{sec:duality-critical-lines}
	
We first show that~$\atrc_{J,U}^{0,0}$ and $\atrc_{J,U}^{1,1}$ coincide for almost every~$(J,U)$:

\begin{lemma}\label{lemma:dense_subsets}
	There exists $D\subseteq\{(J,U)\in\R^2:0<J<U\}$ with Lebesgue measure~$0$ such that, for any $(J,U)\in D^c$, one has
	\begin{equation}\label{eq:free=wired}
	\atrc_{J,U}^{0,0}=\atrc_{J,U}^{1,1}.
	\end{equation}
\end{lemma}

\begin{remark}
Note that, by~\eqref{eq:cbc}, equation~\eqref{eq:free=wired} implies equality of all Gibbs measures. 
\end{remark}

The proof goes by applying the classical FK-percolation argument to the marginals of~$\atrc$ on~$\omega_\tau$ and~$\omega_{\tau\tau'}$, see Appendix~\ref{sec:dense_subsets} for more details.
We are ready to prove part (i) of Theorem~\ref{thm:duality}. Recall that $B_{r}(x,y)$ is the Euclidean ball of radius $r$ centred at $(x,y)$.
\begin{proof}[Proof of Theorem~\ref{thm:duality}(i)]
Fix $J<U$. By~\eqref{eq:cbc}, we have $\betac^{\tau}\leq\betac^{\tau,f}$. Assume for contradiction that the inequality is strict, and take $\beta\in(\betac^{\tau},\betac^{\tau,f})$. Then, by~\eqref{eq:strong-mon}, 
there exists $\eps>0$ such that, for any $(J',U')\in B_\eps(\beta J,\beta U)$, 
\[
	\atrc_{J',U'}^{0,0}(0\xlra{\omega_\tau}\infty)=0\quad\text{and}\quad \atrc_{J',U'}^{1,1}(0\xlra{\omega_\tau}\infty)>0. 
\]
This contradicts Lemma~\ref{lemma:dense_subsets}.
\end{proof}

Denote by $\calH^{\tau}_n$ (resp. $\calH^{\tau\tau'}_n$) the event that the box $[0,2n-1]\times[0,2n-1]$ is crossed horizontally by $\omega^{\tau}$ (resp. $\omega^{\tau\tau'}$). 
Note that the complement of $\calH^{\tau}_n$ is the event that the box $[0,2n-1]\times[0,2n-1]$ is crossed vertically by the dual $\omega^{*}_{\tau}$.
The following lemma states a standard characterisation of non-transition points.
It is a consequence of Lemma~\ref{lemma:dense_subsets} and sharpness of the phase transition in the ATRC.
The latter can be derived using a robust approach going through the OSSS inequality~\cite{DumRaoTas19}; see Appendix~\ref{sec:sharpness} for more details.
 
\begin{lemma}\label{lemma:atrc_cont.}
	Let $0<J<U$. 
	Then, $(J,U)\in \gamma_\tau$ if and only if, for any $\eps>0$, 
	there exist points~$(J_0,U_0)$ and~$(J_1,U_1)$ in~$B_\eps( J, U)$, such that, as~$n\to \infty$,
	\begin{equation}\label{eq:critical-criterion}
		\atrc_{J_0,U_0}^{0,0}[\calH_n^{\tau}]\to 1 
		\quad\text{and}\quad
		\atrc_{J_1,U_1}^{1,1}[\calH_n^{\tau}]\to 0.
	\end{equation}
	The same holds also when $\tau$ is replaced everywhere by $\tau\tau'$.
\end{lemma}

\begin{proof}
	Assume~$(J,U)\in \gamma_\tau$.
	By sharpness, $\atrc_{\beta J, \beta U}^{1,1}[\calH_n^{\tau}]\rightarrow 0$, for any~$\beta \in (0,1)$.
	Also, it is standard that part~$(i)$ of Theorem~\ref{thm:duality} and Zhang's argument imply that
	$\atrc_{\beta J, \beta U}^{0,0}[\calH_n^{\tau}]\rightarrow 1,$ for any~$\beta>1$.
	This gives one direction of the statement.
	
	To show the reverse, assume first that~$\betac^\tau:=\betac^\tau(J,U)>1$ 
	and take~$\beta \in (1,\betac^\tau)$.
	By sharpness, $\atrc_{\beta J,\beta U}^{1,1}[\calH_n^{\tau}]\rightarrow 0$ and, by~\eqref{eq:strong-mon}, the same holds in some neighbourhood of~$(J,U)$.
	The case~$\betac^\tau < 1$ is analogous.
\end{proof}

\begin{proof}[Proof of Theorem~\ref{thm:duality}(ii)]
	Let~$(J,U)\in\gamma_{\tau}$ and~$\eps >0$.
	Since the duality mapping is a continuous involution, we can find~$\delta>0$ such that the image of~$B_\delta(J,U)$ is inside~$B_\eps(J^*,U^*)$.
	By Lemma~\ref{lemma:atrc_cont.}, we get \eqref{eq:critical-criterion} for some~$(J_0,U_0)$ and~$(J_1,U_1)$ in~$B_\delta(J,U)$.
	By duality and symmetry,
	\begin{equation*}
		\atrc_{J_0^{*},U_0^{*}}^{1,1}[\calH_n^{\tau\tau'}]\to 0
		\quad\text{and}\quad 
		\atrc_{J_1^{*},U_1^{*}}^{0,0}[\calH_n^{\tau\tau'}]\to 1\quad\text{as }n\to\infty.
	\end{equation*}	
	Since~$(J_0^{*},U_0^{*})$ and~$(J_1^{*},U_1^{*})$ are in~$B_\eps(J^*,U^*)$, Lemma~\ref{lemma:atrc_cont.} implies that~$(J^*,U^*) \in \gamma_{\tau\tau'}$.
	
	Proving that~$(J^*,U^*) \in \gamma_{\tau\tau'}$ implies~$(J,U)\in\gamma_{\tau}$ is analogous.
\end{proof}

\subsection{$\varphi_{\beta}(S)$ argument: proof of Theorem~\ref{thm:splitting}}
\label{sec:phi-s}

Following~\cite{Sim80,Lie80} (see also~\cite{DumTas15}), for a finite subgraph $S\subset \bbL$ containing $0$, define 
\begin{equation*}
	\varphi_{\beta}(S)=|\partial S|\cdot \atrc^{1,1}_{S,\beta}(0\xlra{\omega_\tau} \partial S).
\end{equation*}

The following lemma states a key property of~$\varphi_{\beta}(S)$: if it is less than~$1$ for some~$S$, then~$\omega_\tau$ exhibits exponential decay of connection probabilities.
This finite criterion allows to use continuity of~$\varphi_{\beta}(S)$ and Proposition~\ref{thm:finite_expo_decay_at_sd} to extend exponential decay of~$\omega_\tau$ beyond~\eqref{eq:self-dual}. Let $\Lambda_{k}$ be the box of size $k$ in $\bbL$, that is $\Lambda_k=\{u\in\bbL:\norm{u}_1\leq 2k\}$.

\begin{lemma}\label{lemma:phi(S)}
	Let $\beta>0$. Assume that $\varphi_{\beta}(S)<1$, for some finite subgraph $S\subset \bbL$ containing $0$.
	Then, there exists $c:=c(\beta,S)>0$ such that 
	\begin{equation*}
		\atrc^{1,1}_{\beta}(0\xlra{\omega_\tau}\partial\Lambda_{n})\leq e^{-cn}.
	\end{equation*}
\end{lemma}

\begin{remark}
	Note that the boundary conditions are {\em free} in~\cite{DumTas15} and {\em wired} in our case.
	The reason is that an analogue of Lemma~\ref{lemma:phi(S)} is proven in~\cite{DumTas15} via a modified Simon--Lieb inequality~\cite{Lie80, Sim80} for the Ising model.
	Such inequalities are not available in our case.
	While Lemma~\ref{lemma:phi(S)} under wired conditions is elementary, proving exponential decay under wired boundary conditions in finite volume (Proposition~\ref{thm:finite_expo_decay_at_sd}) is the subject of Sections~\ref{sec:couplings}-\ref{sec:exp-dec-finite}.
\end{remark}

\begin{proof}[Proof of Lemma~\ref{lemma:phi(S)}]
	Let $S \subset \bbL$ be a finite subgraph containing $0$ such that $\varphi_{\beta}(S)<1$ and let~$k$ be such that $S\subset\Lambda_{k}$. 
	If $0 \xlra{\omega_\tau}\partial\Lambda_{nk}$, then $\partial S \xlra{\omega_\tau}\partial\Lambda_{nk}$ and~$0\xlra{\omega_\tau} \partial S$.
	
	By~\eqref{eq:cbc} and the union bound,
	\begin{align*}
		\atrc^{1,1}_{\beta}(0\xlra{\omega_\tau}\partial\Lambda_{nk})
		&\leq 
		\atrc_{\beta}^{1,1}(0\xlra{\omega_\tau} \partial S \, | \,
		\partial S \xlra{\omega_\tau}\partial\Lambda_{nk}) \cdot
		\atrc^{1,1}_{\beta}(\partial S\xlra{\omega_\tau}\partial\Lambda_{nk})
		\\
		&\leq
		\atrc_{S,\beta}^{1,1}(0\xlra{\omega_\tau} \partial S) \cdot
		\sum_{x\in\partial S} \atrc^{1,1}_{\beta}(x\xlra{\omega_\tau}\partial\Lambda_{nk})\\
		&\leq 
		\atrc_{S,\beta}^{1,1}(0\xlra{\omega_\tau} \partial S) \cdot
		\abs{\partial S}\atrc_{\beta}^{1,1}(0\xlra{\omega_\tau}\partial\Lambda_{(n-1)k})\\
		&= 
		\varphi_{\beta}(S)\atrc^{1,1}_{\beta}(0\xlra{\omega_\tau}\partial\Lambda_{(n-1)k}).
	\end{align*}
	where we also used translation invariance of~$\atrc_\beta^{1,1}$ and that $S\subset\Lambda_{k}$.

	Since $\varphi_{\beta}(S)<1$, we get that $\atrc_{\beta}^{1,1}(0\xlra{\omega_\tau}\Lambda_{nk})$ decays exponentially fast in $n$ by induction. Since for any $m\in\bbN$, there exists $n$ such that $m\in [nk,(n+1)k]$, we get that $\atrc^{1,1}_{\beta}(0\xlra{\omega_\tau}\partial\Lambda_{m})$ decays exponentially fast in $m$.
\end{proof}

We are now ready to derive Theorem~\ref{thm:splitting} from Proposition~\ref{thm:finite_expo_decay_at_sd} and Theorem~\ref{thm:duality}.

\begin{proof}[Proof of Theorem~\ref{thm:splitting}]
	Fix~$J<U$.
	By Proposition~\ref{thm:finite_expo_decay_at_sd}, we can take~$n>1$ such that $\varphi_{\betasd}(\Lambda_n)<1$.
	Since the function $\beta\mapsto\varphi_{\beta}(\Lambda_n)$ is increasing and continuous, there exists~$\varepsilon=\varepsilon(J,U)>0$, such that~$\varphi_{\beta'}(\Lambda_n)<1$, for all~$\beta'< \betasd+\varepsilon$.
	The latter implies exponential decay by Lemma~\ref{lemma:phi(S)} and hence~$\beta_c^{\tau} > \betasd$.
	In other words, all points on $\gamma^\tau$ are strictly above the self-dual curve. 
	Hence their images under the duality mapping~\eqref{eq:duality_1} are strictly below the self-dual curve.
	By Theorem~\ref{thm:duality}, these points are exactly the points of~$\gamma^{\tau\tau'}$ and this finishes the proof.	
\end{proof}
	
\begin{remark}\label{rem:continuity-argument}
Standard arguments similar to the proof of Lemma~\ref{lemma:phi(S)} show that
\[
c_\beta=\lim_{n\to\infty}-\frac{1}{n}\log\atrc_{\Lambda_n,\beta}^{1,1}[0\xlra{\omega_\tau}\partial\Lambda_n]
\]
exists and is right-continuous in $\beta$, which gives another way to argue that the exponential decay from Proposition~\ref{thm:finite_expo_decay_at_sd} extends to an open neighbourhood of the self-dual line~\eqref{eq:self-dual}. 
\end{remark}

\section{Models, couplings and required input}
\label{sec:couplings}

In this section, we introduce the six-vertex model together with its height and spin representations.
We also state couplings of this model with the ATRC model and FK percolation that will be crucial to our arguments.
A combination of these two couplings has been made explicit recently in the work of Peled and the third author~\cite{GlaPel19} and we summarize the results of that work that we will rely on.

\begin{figure}
	\begin{center}
		\includegraphics[scale=0.8]{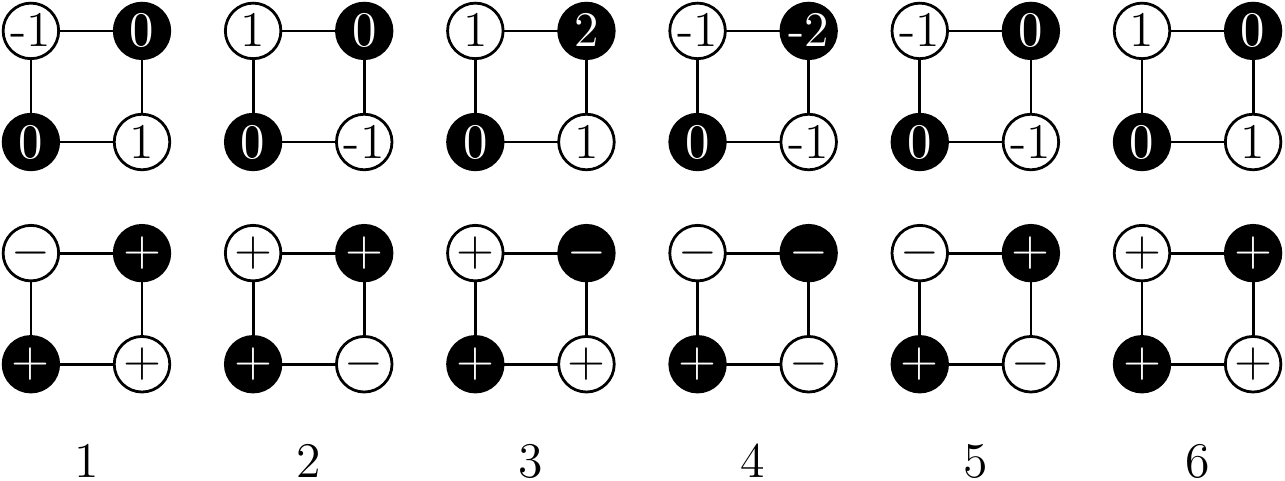}
	\end{center}
	\caption{\emph{Top:} The height representation of the six-vertex model in the four vertices of a unit square in $\bbZ^{2}$, normalized to equal $0$ at the lower left vertex.
\emph{Bottom:} The spin representation is derived from the heights by setting the spin state at each vertex to $+1$ (resp.~$-1$) if the height modulo~$4$ equals $0,1$ (resp.~$2,3$).}
	\label{fig:6v-hom-config}
\end{figure}

\subsection{Graph notation}
\label{sec:notation}

\paragraph{Dual subgraphs and configurations.}
For a finite subgraph~$\Omega$ of~$\bbL$, define its dual graph~$\Omega^*$ in~$\bbL^*$ formed by edges dual to the edges of~$\Omega$. As for primal graphs, we denote the sets of its vertices and edges by~$V_{\Omega^*}$ and~$E_{\Omega^*}$. The boundary~$\partial\Omega^*$ is defined in the same way as for subgraphs of~$\bbL$. Given a percolation configuration~$\omega\in\{0,1\}^{E_\Omega}$, its dual configuration~$\omega^*\in\{0,1\}^{E_{\Omega^*}}$ is defined by~$\omega^*(e^*)=1-\omega(e),\,e\in E_\Omega$.
\paragraph{Domains in $\bbL$.}
A finite subgraph~$\Omega$ of~$\bbL$ (or~$\bbL^*$) is a \emph{domain} if it is induced by vertices within a simple cycle (including the cycle itself). 
We denote the set of vertices on the surrounding cycle by~$\overline{\partial}\Omega$ and call it the \emph{domain-boundary} of~$\Omega$.
The set of edges on~$\overline{\partial}\Omega$ is called \emph{edge-boundary} of~$\Omega$ and is denoted by~$E_{\overline{\partial}\Omega}$.

\paragraph{Domains in $\bbZ^2$.}
Given a domain~$\Omega$ in~$\bbL$, let~$\calD_{\Omega}$ be the subgraph of~$\bbZ^{2}$ induced by vertices in~$\Omega\cup\Omega^{*}$. We call such a domain an \textit{even domain} of~$\bbZ^{2}$ (see Figure~\ref{fig:at-diagram}). 
Define~$\partial\calD_{\Omega}=\overline{\partial}\Omega\cup\partial\Omega^{*}$. 
Given a domain~$\Omega'$ on~$\bbL^{*}$, we define~$\calD_{\Omega'}$ in the same manner and call it an \textit{odd domain} of~$\bbZ^{2}$. 
We emphasize that we only consider even and odd domains. 

\begin{remark}
The reason for the unusual choice of boundary $\partial\calD_{\Omega}=\overline{\partial}\Omega\cup\partial\Omega^{*}$ (rather than $\partial\Omega\cup\partial\Omega^*$) lies in the Baxter--Kelland--Wu coupling (Section~\ref{sec:bkw}).
\end{remark}

\subsection{Six-vertex model and its representations}
\label{sec:sixv}

In this section, we define the \emph{six-vertex} model (more precisely, the \emph{F-model}) and its different representations in terms of
 spins and height functions. For the whole subsection, fix a domain $\Omega$ in $\bbL$ (or $\bbL^*$) and its corresponding even (odd) domain $\calD=\calD_\Omega$ in $\bbZ^2$.

\paragraph{Height functions.}
A function $h:\calD\to\Z$ is called a \emph{height function} (of the six-vertex model) if it satisfies the \emph{ice rule}:
\begin{itemize}
    \item $\abs{h(u)-h(v)}=1$ whenever $u,v$ are connected by an edge in $\bbZ^2$,
    \item $h$ takes even values on~$\calD\cap\bbL$.
\end{itemize}
This constraint implies that, for each edge~$e$, the value of~$h$ is constant either at the endpoints of~$e$ or at the endpoints of~$e^*$.
Up to an even additive constant, this leaves six local possibilities (\emph{types}), where types~$5$ and~$6$ correspond to~$h$ taking constant values along both~$e$ and~$e^*$, see Figure~\ref{fig:6v-hom-config}.

The \emph{six-vertex height function measure} on $\calD$ with parameters $c,c_{b}>0$ and boundary conditions~$t\in\bbZ^{\partial \calD}$ is supported on height functions $h\in\Z^{\calD}$ that coincide with~$t$ on~$\partial\calD$ and is given by
\begin{equation}\label{eq:hf-def}
    \HF_{\calD,c}^{t;c_b}[h]=\frac{1}{Z^{t,c_b}_{\mathsf{hf},\calD,c}} c^{n_{5,6}^{\rm i}(h)} c_b^{n_{5,6}^{\rm b}(h)},
\end{equation}
where $Z^{t,c_b}_{\mathsf{hf},\calD,c}$ is a normalizing constant and $n_{5,6}^{\rm i}(h)$ (resp.~$n_{5,6}^{\rm b}(h)$) is the number of edges of type 5 or 6 of $E_\Omega\setminus E_{\overline{\partial}\Omega}$ (resp.~$E_{\overline{\partial}\Omega}$). When $c=c_{b}$, we recover the standard six-vertex probability measure that will be denoted by $\HF_{\calD,c}^{t}$.
We write~$\HF_{\calD,c}^{2n,2n+1;c_b}$ for~$\HF_{\calD,c}^{t;c_b}$ with $t\in\lbrace 2n,2n+1\rbrace^{\partial \calD}$. 
We define~$\HF_{\calD,c}^{2n,2n-1;c_b}$ analogously.  

Finally, we define $\HF_{\calD,c}^{2n,2n\pm 1;c_b}$ as the probability measure given by~\eqref{eq:hf-def} and supported on all height functions in~$\bbZ^\calD$ that have a fixed value~$2n$ on~$\partial\calD\cap\bbL$.
Note that the value on~$\partial\calD\cap\bbL^*$ is not fixed in this case, so the conditions can be viewed as {\em semi-free}.

\paragraph{Spin representation.}
Given a height function $h\in\Z^{\calD}$, define $\sigma=\sigma(h)\in\{\pm1\}^{\calD}$ by
\[\sigma(u)=
\begin{cases}
   1 & \text{if }h(u)\equiv 0,1\text{ (mod 4)},\\
   -1 & \text{otherwise.}
\end{cases}\]
The \emph{six-vertex spin measures} $\Spin_{\calD,c}^{+,+;c_b}$, $\Spin_{\calD,c}^{+,+}$, $\Spin_{\calD,c}^{+,\pm}$ are defined as the push-forwards of~$\HF_{\calD,c}^{0,1;c_b}$, $\HF_{\calD,c}^{0,1}$, $\HF_{\calD,c}^{0,\pm 1}$ under this mapping.
The spin measures are supported on all spin configurations~$\sigma\in \{\pm 1\}^{\calD}$ with the following restrictions:
$\sigma_{|\partial\calD} \equiv 1$ under~$\Spin_{\calD,c}^{+,+;c_b}$ and $\Spin_{\calD,c}^{+,+}$;
$\sigma_{|\partial\calD\cap \bbL} \equiv 1$ under~$\Spin_{\calD,c}^{+,\pm}$.

\begin{figure}
	\begin{center}
		\includegraphics[scale=0.8, page=1]{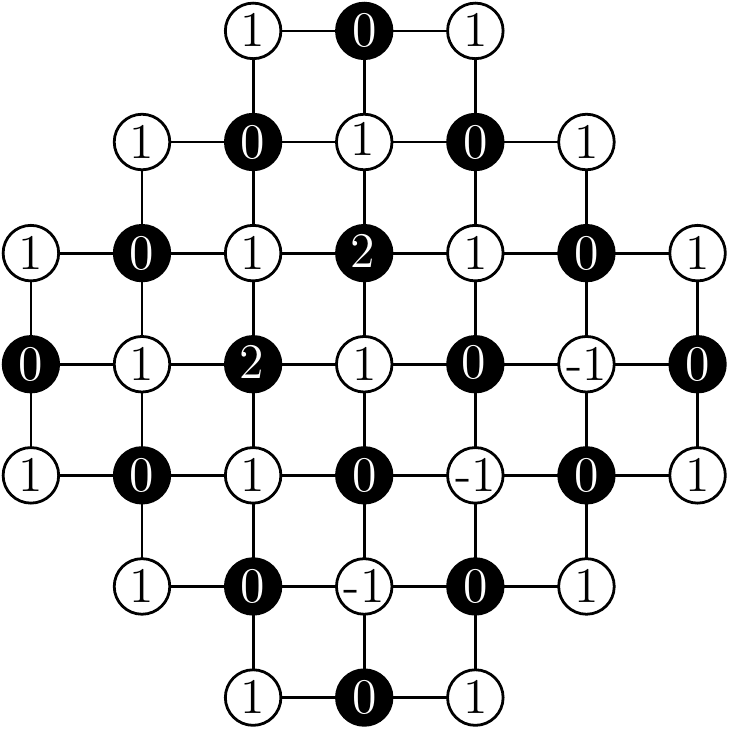}
		\quad\quad\quad
		\includegraphics[scale=0.8, page=2]{heights_spins.pdf}
	\end{center}
	\caption{\emph{Left:} Height function with~$0,1$ boundary conditions. \emph{Right:} Its spin representation is given by~$\sigma^\bullet$ on~$\bbL$ (black circles) and~$\sigma^\circ$ on~$\bbL^*$ (white circles).}
	\label{fig:height-function-0-1}
\end{figure}

\subsection{From Ashkin--Teller to six-vertex}
\label{sec:sixv-at}

In this section, we describe the connection between the self-dual Ashkin--Teller model on a domain $\Omega$ in $\bbL$ and the spin representation of the six-vertex model on the corresponding even domain $\calD_{\Omega}$ in $\bbZ^2$.
We consider two types of boundary conditions that will play an important role in proving Proposition \ref{thm:finite_expo_decay_at_sd}.

Let~$\Omega$ be a domain of~$\bbL$, $J<U$ be parameters. We will consider the ATRC measures defined in Section~\ref{sec:rc-at} with boundary conditions on $\overline{\partial}\Omega$ rather than $\partial\Omega$. We write $\atrc_{\Omega,J,U}^{0,\mathbbm{1}}$ and $\atrc_{\Omega,J,U}^{\mathbbm{1},\mathbbm{1}}$ for the corresponding ATRC measures where $\mathbbm{1}$ refers to the wired boundary condition on $\overline{\partial}\Omega$.

Let~$\eta_\tau,\eta_{\tau\tau'}$ be boundary conditions on $\partial\Omega$ or $\overline{\partial}\Omega$. Consider the marginal of $\atrc_{\Omega,J,U}^{\eta_\tau,\eta_{\tau\tau'}}$ on~$\omega_\tau$: 
this is the probability measure supported on~$\{0,1\}^{E_\Omega}$ and defined by
\[
	\mu_{\Omega,J,U}^{\eta_\tau,\eta_{\tau\tau'}}(\xi):= 
	\atrc_{\Omega,J,U}^{\eta_\tau,\eta_{\tau\tau'}}(\{\omega_\tau = \xi\}).
\]
Also, given $\sigma\in\{\pm 1\}^{\calD_{\Omega}}$, we write $\sigma^\bullet$ and $\sigma^\circ$ for the restrictions of $\sigma$ to $\Omega$ and $\Omega^{*}$, respectively (see Fig~\ref{fig:height-function-0-1}). We define the sets of disagreement edges:
\[
	E_{\sigma^\bullet}:=\{uv\in E_\Omega:\sigma^\bullet(u)\neq \sigma^\bullet(v)\}
	\quad\text{and}\quad
	E_{\sigma^\circ}:=\{u^*v^*\in E_{\Omega^*}:\sigma^\circ(u^*)\neq \sigma^\circ(v^*)\}.
\]
Finally, we define the {\em compatibility} relation on pairs of~$\sigma^\bullet\in\{\pm1\}^{\Omega}$ 
and~$\omega\in\{0,1\}^{E_\Omega}$:
\[ 
	\sigma\perp\omega \quad
	\text{if and only if} \quad
	\sigma(u) = \sigma(v), \text{ for any } uv\in \omega.
\]
The compatibility relation on pairs of~$\sigma^\circ\in\{\pm1\}^{\Omega^*}$ and~$\omega\in\{0,1\}^{E_{\Omega^*}}$ is defined similarly.
The following is a consequence of \cite[Proposition 8.1]{GlaPel19} and a remark after it, or may be proved along the same lines:

\begin{proposition}\label{prop:6v-atrc-coupling}
	Let~$0<J<U$ be a point on the self-dual line~\eqref{eq:self-dual} and~$c=\coth{2h}$.
	
	1) If $\Omega$ is a domain in $\bbL$, then we can couple $\sigma\sim \Spin_{\calD_{\Omega},c}^{+,+}$ and $\omega_\tau\sim\mu_{\Omega,J,U}^{0,\mathbbm{1}}$ by
	   	\begin{align*}
	    \Prob{\sigma,\omega_\tau}\propto\left(\tfrac{1}{c-1}\right)^{\abs{\omega_\tau}+\abs{E_{\sigma^\bullet}}}\ind{\sigma^\bullet\perp\omega_\tau,\,\sigma^\circ\perp\omega_\tau^*}\ind{\sigma^\bullet\equiv +\text{ on }\overline{\partial}\Omega,\,\sigma^\circ\equiv +\text{ on }\partial\Omega^{*}}.
	    \end{align*}
	    Thus, $\sigma^\circ$ is obtained by assigning $+1$ to the clusters of~$\omega_\tau^*$ that intersect~$\partial\Omega^*$ and assigning $\pm 1$ uniformly independently to all other clusters.
	    
	2) If $\Omega^*$ is a domain in $\bbL^*$, then we can couple $\sigma\sim \Spin_{\calD_{\Omega^*},c}^{+,\pm}$ and $\omega_\tau\sim\mu_{\Omega,J,U}^{\mathbbm{1},\mathbbm{1}}$ by
	    \begin{align*}
	    \Prob{\sigma,\omega_\tau}\propto\left(\tfrac{1}{c-1}\right)^{\abs{\omega_\tau}+\abs{E_{\sigma^\bullet}}}\ind{\sigma^\bullet\perp\omega_\tau,\,\sigma^\circ\perp\omega_\tau^*}\ind{\sigma^\bullet\equiv +\text{ on }\partial\Omega}.
	    \end{align*}
\end{proposition}

\begin{remark}
	Part 1) of Proposition~\ref{prop:6v-atrc-coupling} is a special case of \cite[Proposition 8.1]{GlaPel19} while part 2) may be proved in the same way. The proof relies on the following identity:
	\[
	k^1(\omega)-k(\omega^*)=\abs{\omega^*}+\mathrm{const}(\Omega),
	\]
	where $\omega=\omega_\tau^*$ for 1) and $\omega=\omega_\tau$ for 2).
	This follows from Euler's formula using that either $\Omega$ or $\Omega^*$ is a domain and our definition of the domain-boundary.
\end{remark}

\begin{corollary}\label{cor:omega-tau-from-spins}
	 In the setting of part 1) of Proposition~\ref{prop:6v-atrc-coupling}, take $(\sigma^\bullet,\sigma^\circ)\sim\Spin_{\calD_{\Omega},c}^{+,+}$.
	 Sample a percolation configuration~$\omega$ on~$E_\Omega$ as follows independently at each edge~$e$: if the endpoints of~$e^*$ have opposite values in~$\sigma^\circ$, then $\omega_e = 1$; if the endpoints of~$e$ have opposite values in~$\sigma^\bullet$, then $\omega_e = 0$; if~$\sigma^\circ$ agrees on~$e^*$ and~$\sigma^\bullet$ agrees on~$e$, then
	\begin{equation}\label{eq:ber-c}
		\prob(\omega_e = 1) = \tfrac1{c}.
	\end{equation}
	Then, the law of~$\omega$ is given by~$\mu_{\Omega,J,U}^{0,\mathbbm{1}}$.
\end{corollary}

\subsection{Input from the six-vertex model}
\label{sec:sixv-input}

In this section, we mention basic properties of six-vertex measures and state some results from \cite{GlaPel19}. The following proposition is a combination of~Theorem~2, Proposition~6.1 and Lemma~6.2 in~\cite{GlaPel19}. 
We remark that we only consider even and odd domains in $\bbZ^2$.

For~$u\in\bbL$, $S\subset \bbL$, define~$u\xleftrightarrow{h\neq 0}S$ to be the event that $u$ is connected (in~$\bbL$) to $S$ by a path of heights different from~$0$.
We similarly define~$u^* \xleftrightarrow{h\neq 1}S^*$ for~$u^*\in\bbL^*$, $S^*\subset \bbL^*$.

\begin{proposition}[\cite{GlaPel19}]\label{proposition:HF_weak_limit}
	Fix $c>2$, and let $\lambda$ be the unique positive solution of $c=e^{\lambda/2}+e^{-\lambda/2}$. 
	Then, for any sequence of domains $\calD_{k}\nearrow\bbZ^{2}$, the measures $\HF_{\calD_{k},c}^{0,1}$ and $\HF_{\calD_{k},c}^{0,1;e^{\lambda/2}}$ converge weakly to the same limit that we denote by $\HF_{c}^{0,1}$.
	Moreover, $\HF_{c}^{0,1}$-a.s. exist unique infinite clusters in~$\bbL$ of height~$0$ and in~$\bbL^*$ of height~$1$.
	Finally, clusters of other heights are exponentially small: for some $\alpha>0$ uniform in~$n\geq 1$, $u\in \bbL$, $u^*\in\bbL^*$,
	\[
		\HF_{c}^{0,1}\left[u\xleftrightarrow{h\neq 0}u+\partial\Lambda_n\right]<e^{-\alpha n} \quad \text{and} \quad
		\HF_{c}^{0,1}\left[u^*\xleftrightarrow{h\neq 1}u^*+\partial\Lambda_n\right]<e^{-\alpha n}.
	\]
\end{proposition}

Let us emphasize that, while existence of subsequential limits is a straightforward consequence of discontinuity of the phase transition in FK-percolation, the ordering of {\em both} even and odd heights is non-trivial.
This also implies that the weak limit of~$\HF_{\calD_k,c}^{0,1;e^{\lambda/2}}$ remains the same, whether it is taken along even or odd domains.
Analogously, for any $n\in\bbZ$, one obtains limit measures $\HF_{c}^{2n,2n+1}$ (resp. $\HF_{c}^{2n,2n-1}$) of $\HF_{\calD_k,c}^{2n,2n+1}$ (resp. $\HF_{\calD_k,c}^{2n,2n-1}$) satisfying the corresponding properties. 

Since the modulo 4 mapping (Section~\ref{sec:sixv}) is local, Propositon~\ref{proposition:HF_weak_limit} directly implies the following corollary.

\begin{corollary}\label{cor:spin-weak-limits}
	Fix $c>2$, and let $\lambda$ be the unique positive solution of $c=e^{\lambda/2}+e^{-\lambda/2}$. Then, for any sequence of domains $\calD_{k}$ increasing to $\bbZ^{2}$, the measures $\Spin_{\calD_{k},c}^{+,+}$ and $\Spin_{\calD_{k},c}^{+,+;e^{\lambda/2}}$ converge weakly to some $\Spin_{c}^{+,+}$, which is independent of the sequence $\calD_k$.
\end{corollary}

The height function measures admit useful monotonicity properties and correlation inequalities when $c,c_b\geq 1$, see \cite[Proposition 5.1]{GlaPel19}.

\begin{proposition}\label{prop:hf-fkg}
	Let $\calD$ be a domain in $\bbZ^{2}$, and let $c,c_b\geq 1$. Then, for any boundary condition $t$, the measure $\HF_{\calD,c}^{t;c_b}$ satisfies the FKG inequality~\eqref{eq:fkg-general}.
	In particular, if $t\leq t'$, then $\HF_{\calD,c}^{t;c_b}$ is stochastically dominated by $\HF_{\calD,c}^{t';c_b}$.
\end{proposition}

It has been established in~\cite[Theorem 4]{GlaPel19} that the marginals of~$\Spin_{\calD_{k},c}^{+,+}$ on~$\sigma^\bullet$ (resp.~$\sigma^\circ$) satisfy the FKG inequality with respect to the pointwise order on~$\{\pm1\}^{V(\Omega)}$ (resp.~$\{\pm1\}^{V(\Omega^*)}$).
Though~\cite{GlaPel19} deals only with boundary conditions specified on the whole boundary, the extension to free or semi-free conditions is straightforward.
Indeed, the statement for~$\sigma^\bullet$ holds as long as spins~$\sigma^\circ$ on the boundary are not forced to disagree.

\begin{proposition}[\cite{GlaPel19}]\label{prop:spin-fkg}
	Let $\Omega$ be a domain in $\bbL$, and let $c\geq 1$. 
	The marginals of $\Spin_{\calD_{\Omega},c}^{+,\pm}$ on $\sigma^\bullet$ and~$\sigma^\circ$ satisfy the FKG inequality~\eqref{eq:fkg-general}.
\end{proposition}

It was also shown in \cite[Corollary 7.3, Proposition 7.5]{GlaPel19} that the marginals $\mu_{\Omega,J,U}^{0,\mathbbm{1}}$ converge to some infinite volume state $\mu_{J,U}^{0,1}$ that admits exponential decay of connection probabilities.

\begin{proposition}[\cite{GlaPel19}]\label{prop:exp.decay_atrc01}
	Let $0<J<U$ be on the self-dual line~\eqref{eq:self-dual} and $\Omega_k$ be a sequence of domains increasing to $\bbL$. The measures $\mu_{\Omega_{k},J,U}^{0,\mathbbm{1}}$ converge weakly to some measure $\mu_{J,U}^{0,1}$ on $\{0,1\}^{E(\bbL)}$ which is independent of the sequence $\Omega_k$ and admits exponential decay of $\omega_\tau$-connection probabilities: there exist $M,\alpha>0$ such that, for any $u,v\in\bbL$,
	\[
		\mu_{J,U}^{0,1}[u\xleftrightarrow{\omega_\tau}v]\leq Me^{-\alpha\abs{u-v}}.
	\]
\end{proposition}
We sketch the argument given in~\cite{GlaPel19}.

\begin{proof}[Sketch of proof of Proposition \ref{prop:exp.decay_atrc01}]
	The couplings in Proposition \ref{prop:6v-atrc-coupling} and Corollary \ref{cor:spin_conv} imply convergence of $\mu_{\Omega,J,U}^{0,\mathbbm{1}}$, as $\Omega\nearrow\bbL$, to some $\mu_{J,U}^{0,1}$ that satisfies FKG and is invariant to translations.
	Thus, it is enough show that that it is exponentially unlikely that~$\omega_\tau$ contains a circuit surrounding~$\Lambda_n$.
	Indeed, on this event, the marginal of $\Spin_c^{+,+}$ at vertices in~($\Lambda_n)^*$ is invariant to the spin flip.
	By Proposition \ref{exp.relax.HF_even}, radii of clusters of minuses have exponential tails and the claim follows.
\end{proof}

We emphasise a difference between Propositions~\ref{prop:exp.decay_atrc01} and~\ref{thm:finite_expo_decay_at_sd}: the latter proves exponential decay under the largest boundary conditions and in finite volume.
As we saw in Section~\ref{sec:phi-s}, this is necessary for the proof of Theorem~\ref{thm:splitting}.

\subsection{FK-percolation}
\label{sec:fk}

Fortuin--Kasteleyn (FK) percolation~\cite{ForKas72} is an archetypical dependent percolation model. It is well-understood thanks to recent remarkable works; see~\cite{Dum17a,Gri06} for background. We will transfer some known results from FK-percolation to the six-vertex model via the BKW coupling (Section~\ref{sec:bkw}) and further to the self-dual ATRC via the coupling in Proposition~\ref{prop:6v-atrc-coupling}. 

\begin{definition}
	Let~$\Omega\subset\bbL$ be a finite subgraph and $\xi$ a partition of $\partial\Omega$. FK-percolation on~$\Omega$ with parameters~$p\in [0,1]$ and $q>0$ is supported on percolation configurations $\eta\in\{0,1\}^{E_\Omega}$ and is given by
	\[
		\FK_{\Omega,p,q}^{\xi}(\eta)=\frac{1}{Z} p^{\abs{\eta}}(1-p)^{\abs{E_\Omega}-\abs{\eta}}q^{k^{\xi}(\eta)},
	\]
where $Z=Z(\Omega,p,q,\xi)$ is a normalizing constant and $k^{\xi}(\eta)$ was defined in Section~\ref{sec:rc-at}.
\end{definition}
The {\em free} and {\em wired} FK-percolation measures~$\FK_{\Omega,p,q}^{\rm f}$ 
and~$\FK_{\Omega,p,q}^{\rm w}$ are defined by free and wired boundary conditions, respectively (as in Section~\ref{sec:rc-at}). 

We now review several fundamental results about FK-percolation.

\begin{proposition}\label{prop:rc_basic_properties}
	Let $p\in[0,1]$, $q>1$ and~$\Omega_{k}\nearrow \bbL$ be a sequence of subgraphs.
	Then, the weak limits of~$\FK_{\Omega_{k},p,q}^\mathrm{f}$ and~$\FK_{\Omega_{k},p,q}^\mathrm{w}$ exist and do not depend on the chosen sequence:
	\[
		\FK_{p,q}^\mathrm{f}:= \lim_{k\to \infty} \FK_{\Omega_{k},p,q}^\mathrm{f}
		\hspace{5mm} \text{and} \hspace{5mm}
		\FK_{p,q}^\mathrm{w}:= \lim_{k\to \infty} \FK_{\Omega_{k},p,q}^\mathrm{w}.
	\]
	Moreover, these measures are extremal, invariant to translations and satisfy the following ordering, for any finite subgraph $\Omega\subset\bbL$,
	\[
		\FK_{\Omega,p,q}^\mathrm{f}\leq_{\mathrm{st}}\FK_{p,q}^\mathrm{f}\leq_{\mathrm{st}}\FK_{p,q}^\mathrm{w}\leq_{\mathrm{st}}\FK_{\Omega,p,q}^\mathrm{w}. 
	\]
\end{proposition}

As we will see below, the self-dual AT model with~$J<U$ corresponds to FK-percolation with~$q>4$ at~$p=p_{\rm sd}$, where
\[
	p_{\rm sd}:= \tfrac{\sqrt{q}}{\sqrt{q} + 1}.
\]
This model is self-dual: if~$\omega$ has law~$\FK_{p_{\rm sd},q}^{\rm f}$, then~$\omega^*(e^*):=1-\omega(e)$ has law~$\FK_{p_{\rm sd},q}^{\rm w}$.

\begin{theorem}[\cite{DumGagHar16}]\label{Thm:critical_rc_properties}
	Let $q>4$. Then, $\FK_{p_{\rm sd},q}^\mathrm{f}\neq\FK_{p_{\rm sd},q}^\mathrm{w}$ and
\begin{enumerate}[label=(\roman*)]
    \item $\FK_{p_{\rm sd},q}^\mathrm{w}(\text{there exists a unique infinite cluster})=1$,
    \item there exists $\alpha>0$ such that $\FK_{p_{\rm sd},q}^\mathrm{f}(0\leftrightarrow\partial \Delta_n)\leq e^{-\alpha n}$, for any $n\geq 1$, .
\end{enumerate}
\end{theorem}

The second item of this theorem implies \emph{exponential relaxation} at~$p_{\rm sd}$:

\begin{lemma}\label{lemma:exp.relax.rc_pc}
	Let $q>4$. Then, there exists $\alpha>0$ such that, for $n\geq 1$ and any finite subgraph $\Omega\subset\bbL$ that contains~$\Lambda_{2n}$,
	\[
		d_{\mathrm{TV}}\left(\FK_{\Omega,p_{\rm sd}, q}^\mathrm{w}\vert_{\Lambda_n},\FK_{p_{\rm sd},q}^\mathrm{w}\vert_{\Lambda_n}\right)<e^{-\alpha n}.
	\]
\end{lemma}
The proof is standard and goes through the monotone coupling; see Appendix~\ref{sec:exp-rel-fk}.

\subsection{Baxter--Kelland--Wu (BKW) coupling}
\label{sec:bkw}

FK-percolation and the six-vertex model were related to each other for the first time by Temperley and Lieb~\cite{LieTem71} on the level of partition functions.
BKW~\cite{BaxKelWu76} turned this relation into a probabilistic coupling when~$c>2$.
We follow~\cite{GlaPel19} and describe this coupling using a modified boundary coupling constant~$c_b$.

Take~$q>4$, $p=p_{\rm sd}$.
Let~$\lambda >0$ be the unique positive solution to
\[
	e^\lambda + e^{-\lambda} = \sqrt{q}, \text{ and set } c:=e^{\lambda/2} + e^{-\lambda/2}.
\]
Let $\Omega$ be a domain in $\bbL$ and recall the notations introduced in Section~\ref{sec:notation}. The measure $\FK_{\Omega,p_{\rm sd},q}^{\overline{\mathrm{w}}}$ refers to the FK measure with wired boundary conditions on $\overline{\partial}\Omega$. Note that the statements of Proposition~\ref{prop:rc_basic_properties} and Lemma~\ref{lemma:exp.relax.rc_pc} remain valid if we replace~$\mathrm{w}$ by~$\overline{\mathrm{w}}$.

\begin{figure}
	\begin{center}
		\includegraphics[scale=0.8,page=1]{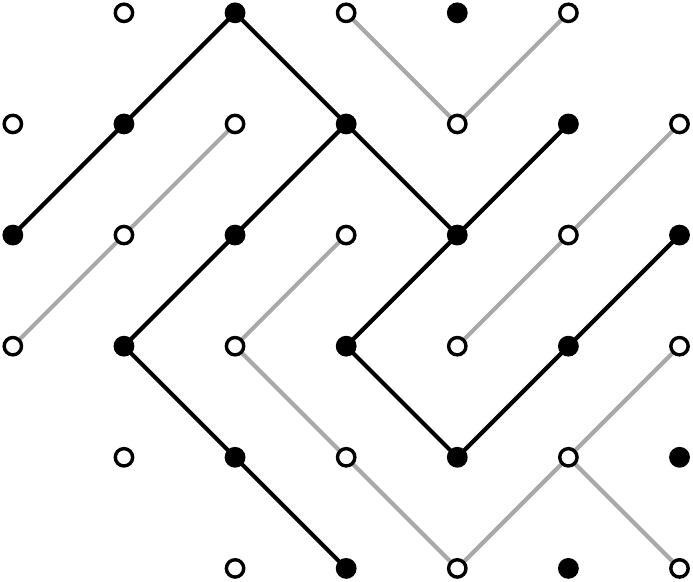}
		\quad\quad\quad
		\includegraphics[scale=0.8,page=2]{BKW.pdf}
	\end{center}
	\caption{\emph{Left:} An edge configuration on the domain $\Omega\subseteq\bbL$ from Figure~\ref{fig:at-diagram} (in black), and its dual on $\Omega^*$ (in gray). \emph{Right:} Loops (in red) separating primal and dual clusters within $\calD_\Omega$ after opening all edges in $E_{\overline{\partial}\Omega}$ (dashed).}
	\label{fig:bkw_loops}
\end{figure}

Consider $\eta\sim\FK_{\Omega,p_{\rm sd},q}^{\overline{\mathrm{w}}}$ and draw loops separating primal and dual clusters within $\Omega$ as in Figure~\ref{fig:bkw_loops}. Given this loop configuration, we define a height function $h\in\Z^{\calD_{\Omega}}$ by:
\begin{itemize}
    \item[\textbf{H1}] Set $h=0$ on $\partial\calD_\Omega\cap\bbL$ and $h=1$ on $\partial\calD_\Omega\cap\bbL^{*}$;
    \item[\textbf{H2}] Assign constant heights to primal and dual clusters by going from~$\partial\calD_{\Omega}$ inside of~$\calD_{\Omega}$ and tossing a coin when crossing a loop: the height {\bf increases} by $1$ with probability $e^{\lambda}/\sqrt{q}$ and {\bf decreases} by $1$ with probability $e^{-\lambda}/\sqrt{q}$, independently of one another.
\end{itemize}

The following result is classical; see e.g. \cite[Chapter 3]{GlaPel19} for a proof in this setup.
\begin{proposition}[BKW coupling]
The resulting height function is distributed according to $\HF_{\calD_{\Omega},c}^{0,1;e^{\lambda/2}}$.
\end{proposition}

\paragraph{Odd domains.} Note that, by symmetry, the whole procedure also works on odd domains with the difference that one needs to replace \textbf{H2} by
\begin{itemize}   
    \item[\textbf{H2'}] each time one crosses a loop, the height \textbf{decreases} by $1$ with probability $e^{\lambda}/\sqrt{q}$ and \textbf{increases} by $1$ with probability $e^{-\lambda}/\sqrt{q}$, independently of one another.
\end{itemize}

For a domain~$\Omega'$ in $\bbL^{*}$, this gives a coupling of $\FK_{\Omega',p_{\rm sd},q}^{\overline{\mathrm{w}}}$ and $\HF_{\calD_{\Omega'},c}^{0,1;e^{\lambda/2}}$.

\section{Exponential decay for ATRC in finite volume}
\label{sec:exp-rel-0-1}

The goal of this section is to derive exponential decay of connection probabilities for $\mu_{\Omega,J,U}^{0,1}$, which is the marginal of a {\em finite-volume} ATRC measure on~$\omega_\tau$, see Section~\ref{sec:sixv-at}. 

Recall that $\Lambda_n=\{u\in\bbL:\norm{u}_1\leq 2n\}$ is the box of size $n$ in $\bbL$.

\begin{proposition}\label{prop:exp-decay-0-1}
	Let $0<J<U$ be on the self-dual line~\eqref{eq:self-dual}.
	There exists $\alpha>0$ such that, for any~$n\geq 1$ and any domain $\Omega$ in $\bbL$ containing $\Lambda_{4n}$,
	\[
		\mu_{\Omega,J,U}^{0,1}[0\xleftrightarrow{\omega_\tau}\partial\Lambda_n]\leq e^{-\alpha n}.
	\]  
\end{proposition}

The proof consists of several steps.
We first transfer the exponential relaxation property from FK-percolation (Lemma~\ref{lemma:exp.relax.rc_pc}) to the six-vertex height function (Proposition~\ref{exp.relax.HF_even}) and then to the marginal~$\nu_{\Omega}$ of the ATRC model with modified edge weights on the boundary. Using Proposition~\ref{proposition:HF_weak_limit}, we also show that the limit of~$\nu_{\Omega}$ is given by~$\mu_{J,U}^{0,1}$ (Lemma~\ref{lemma:nu-convergence}), and that $\nu_{\Omega}$ dominates $\mu_{\Omega,J,U}^{0,1}$ (Lemma~\ref{lemma:stochastic_domination_nu}). The statement then follows from exponential decay in~$\mu_{J,U}^{0,1}$ (Proposition~\ref{prop:exp.decay_atrc01}).

\subsection{Exponential relaxation for height function measures}
\label{sec:exp-rel-hf}

As in Section \ref{sec:bkw}, fix $c>2$ and let $\lambda>0$ be the unique positive solution of $c=e^{\lambda/2}+e^{-\lambda/2}$. Set $q:=(e^\lambda+e^{-\lambda})^2$ and consider $p = p_{\rm sd}$. For $n\geq 1$, define an even domain $\Delta_{2n}=\calD_{\Lambda_n}$.

\begin{proposition}\label{exp.relax.HF_even}
    The convergence of $\HF_{\calD,c}^{0,1;e^{\lambda/2}}$ towards $\HF_{c}^{0,1}$ admits exponential relaxation: there exists $\alpha>0$ such that, for any $n\geq 1$ and any even domain $\calD \supset \Delta_{8n}$,
    \begin{equation}
    d_{\mathrm{TV}}\left(\HF_{\calD,c}^{0,1;e^{\lambda/2}}\vert _{\Delta_{2n}},\HF_{c}^{0,1}\vert_{\Delta_{2n}}\right)<e^{-\alpha n}.\label{eq:exp.relax.HF_even}
    \end{equation}
\end{proposition}

\begin{proof}
	We omit~$q,	p_{\rm sd}$ from the notation for brevity.
	We first construct the limiting measure~$\HF_{c}^{0,1}$. 
	Consider $\eta\sim\FK^{\rm w}$ on $\bbL$.
	Using known results about $\eta\sim\FK^{\rm w}$ (Section~\ref{sec:fk}) we can sample a height function $h$ as follows. 
	Set $h=0$ on the unique infinite cluster of $\eta$ and sample $h$ in its holes according to \textbf{H2} in the BKW coupling (Section~\ref{sec:bkw}). 
	
	Define $\Cc_n$ to be the outermost circuit in~$\eta$ surrounding $\Lambda_{n}$ and contained in $\Lambda_{2n}$ (if it does not exist, we set~$\Cc_n:=\emptyset$).
	Exponential decay in $\eta^*$ and stochastic ordering of FK measures imply existence of $\alpha'>0$ such that, for any $n\geq 1$ and any domain $\Omega \supset \Lambda_{2n}$,
	\begin{equation}
	    \FK_{\Omega}^{\overline{\mathrm{w}}}[\Cc_n\neq\varnothing,\,\Cc_n\xleftrightarrow{\eta}\partial\Omega]\geq \FK^\mathrm{w}[\Cc_n\neq\varnothing,\,\Cc_n\xleftrightarrow{\eta}\infty]>1-e^{-\alpha' n}.\label{RC estimates}
	\end{equation}
	By exponential relaxation of the wired FK measures (Lemma \ref{lemma:exp.relax.rc_pc}), there exists $\alpha''>0$ such that, for any $n\geq 1$ and any domain $\Omega\supset\Lambda_{4n}$,
	\begin{equation}\label{eq:total-var-fk-wired}
		d_{\mathrm{TV}}\left(\FK_{\Omega}^{\overline{\mathrm{w}}}\vert_{\Lambda_{2n}},\FK^\mathrm{w}\vert_{\Lambda_{2n}}\right)<e^{-\alpha'' n}.
	\end{equation}
	Now, given $\Cc_n=C$ and $C\xleftrightarrow{\eta}\infty$, the law of $h$ within $C$ is precisely $\HF_{\calD_{\Omega(C)},c}^{0,1;e^{\lambda/2}}$ where $\Omega(C)$ is the domain in $\bbL$ induced by the vertices within $C$ (including $C$). Note that $\Omega(C)$ contains $\Lambda_n$, whence $\calD_{\Omega(C)}$ contains $\calD_{\Lambda_n}=\Delta_{2n}$.
	
	We can also obtain $\HF_{\calD_{\Omega(C)},c}^{0,1;e^{\lambda/2}}$ from $\FK_{\Omega}^{\overline{\rm w}}$ conditioned on $\Cc_n=C$ and $C\xleftrightarrow{\eta}\partial\Omega$ by applying \textbf{H1} and \textbf{H2}. 
	Together with~\eqref{eq:total-var-fk-wired} and~\eqref{RC estimates}, this proves exponential relaxation.
\end{proof}

Recall that the six-vertex spin measures introduced in Sections \ref{sec:sixv} and \ref{sec:sixv-input} are the push-forwards of the height function measures under the local modulo 4 mapping.

\begin{corollary}\label{cor:spin_conv}
	The convergence of $\Spin_{\calD,c}^{+,+;e^{\lambda/2}}$ towards $\Spin_{c}^{+,+}$ admits exponential relaxation: there exists $\alpha >0$ such that, for any $n\geq 1$ and any even domain $\calD\supset \Delta_{8n}$,
	\begin{equation}
		d_{\mathrm{TV}}\left(\Spin_{\calD,c}^{+,+;e^{\lambda/2}}\vert _{\Delta_{2n}},\Spin_{c}^{+,+}\vert_{\Delta_{2n}}\right)\leq e^{-\alpha n}.\label{eq:exp.relax._spin}
	\end{equation}
\end{corollary}

\subsection{A modified ATRC marginal}\label{sec:modified-atrc}

Fix~$J < U$ on the self-dual line~\eqref{eq:self-dual}, take $c=\coth 2J$ and the unique $\lambda>0$ such that $c=e^{\lambda/2}+e^{-\lambda/2}$.
Let~$\Omega$ be a domain on~$\bbL$ and $\calD_\Omega$ be the corresponding even domain on~$\bbZ^2$.
Recall the definition of the edge-boundary~$E_{\overline{\partial}\Omega}$ in Section~\ref{sec:notation}. Sample $(\sigma^\bullet,\sigma^\circ)$ from~$\Spin_{\calD_{\Omega},c}^{+,+;e^{\lambda/2}}$.
Define~$\nu_\Omega$ as the distribution of $\omega\in\{0,1\}^{E_\Omega}$ sampled independently for each edge~$e$ as follows: if the endpoints of~$e^*$ have opposite values in~$\sigma^\circ$, then $\omega_e = 1$; if the endpoints of~$e$ have opposite values in~$\sigma^\bullet$, then $\omega_e = 0$; if~$\sigma^\circ$ agrees on~$e^*$ and~$\sigma^\bullet$ agrees on~$e$, then
\begin{equation}\label{eq:ber-c-b}
	\prob(\omega_e = 1) =
	\begin{cases}
		e^{-\lambda/2}, &\text{ if } e\in E_{\overline{\partial}\Omega},\\
		\tfrac1{c}, &\text{ otherwise}.
	\end{cases}
\end{equation}

We call~$\nu_\Omega$ a {\em modified ATRC marginal} as it converges to~$\mu_{J,U}^{0,1}$ as~$\Omega\nearrow \bbL$. Moreover, this convergence admits exponential relaxation, which is the content of the next lemma.

\begin{lemma}\label{lemma:nu-convergence}
	For any sequence of domains $\Omega_{k}$ increasing to $\bbL$, the measures $\nu_{\Omega_{k}}$ converge to $\mu_{J,U}^{0,1}$. Moreover, this convergence admits exponential relaxation: there exists $\alpha>0$ such that, for any $n\geq 1$ and any domain $\Omega\supset\Lambda_{4n}$,
	\begin{equation}
		d_{\mathrm{TV}}\left(\nu_{\Omega}\vert _{\Lambda_n},\mu^{0,1}_{J,U}\vert_{\Lambda_n}\right)<e^{-\alpha n}.\label{eq:exp.relax._nu}
	\end{equation}		
\end{lemma}

Recall the representation~\eqref{eq:atrc-repr.} of the ATRC measures. The previous lemma becomes more clear once we identify $\nu_{\Omega}$ as the marginal of $\atrc_{{\rm w}_\tau,{\rm w}_{\tau\tau'}}^{0,\mathbbm{1}}$ on $\omega_\tau$ where the weights are as in~\eqref{eq:atrc-weights-2} except that ${\rm w}_\tau$ is modified on the edge-boundary $E_{\overline{\partial}\Omega}$.

\begin{lemma}\label{lemma:nu-atrc}
	For any domain $\Omega$ in $\bbL$, the measure $\nu_\Omega$ coincides with the marginal of $\atrc_{\Omega,{\rm w}_\tau,{\rm w}_{\tau\tau'}}^{0,\mathbbm{1}}$ on $\omega_\tau$, where  

\begin{equation}\label{eq:nu-atrc-weights}
	\begin{gathered}
		{\rm w}_\tau(e)=
			\begin{cases}
				2\tfrac{c-1}{e^{\lambda/2}-1}, &\text{ if } e\in E_{\overline{\partial}\Omega},\\
				2, &\text{ otherwise},			
			\end{cases}
		\quad\text{and}\quad {\rm w}_{\tau\tau'}\equiv e^{2(U-J)}-1\text{ on }E_\Omega.
	\end{gathered}
\end{equation}
\end{lemma}

We now derive Proposition~\ref{prop:exp-decay-0-1} from Lemmata~\ref{lemma:nu-convergence} and~\ref{lemma:nu-atrc}.
First of all, the ATRC measure is stochastically increasing in ${\rm w}_\tau$ (see Appendix~\ref{sec:proof-stoch-dom-atrc} for the proof).

\begin{lemma}\label{lemma:stoch.dom.atrc}
	Let $\Omega$ be a domain in $\bbL$. The measures $\atrc_{\Omega,{\rm w}_\tau,{\rm w}_{\tau\tau'}}^{0,\mathbbm{1}}$ are stochastically increasing in ${\rm w}_{\tau}$. More precisely, if ${\rm w}_\tau(e)\leq \widetilde{{\rm w}_\tau}(e)$ for all $e\in E_\Omega$, then the measure $\atrc_{\Omega,{\rm w}_\tau,{\rm w}_{\tau\tau'}}^{0,\mathbbm{1}}$ is stochastically dominated by $\atrc_{\Omega,\widetilde{{\rm w}_\tau},{\rm w}_{\tau\tau'}}^{0,\mathbbm{1}}$.
\end{lemma}

This, together with Lemma~\ref{lemma:nu-atrc}, implies the following stochastic domination:

\begin{lemma}\label{lemma:stochastic_domination_nu}
	For any domain $\Omega$ in $\bbL$, the measure $\nu_{\Omega}$ stochastically dominates $\mu_{\Omega,J,U}^{0,\mathbbm{1}}$.
\end{lemma}

\begin{proof}[Proof of Proposition~\ref{prop:exp-decay-0-1}]
	Fix $n\geq 1$ and a domain $\Omega \supset \Lambda_{4n}$ in $\bbL$.
	We have
	\[
		\mu_{\Omega,J,U}^{0,1}[0\leftrightarrow\partial\Lambda_n] \leq \mu_{\Omega,J,U}^{0,\mathbbm{1}}[0\leftrightarrow\partial\Lambda_n] \leq \nu_{\Omega}[0\leftrightarrow\partial\Lambda_n],
	\]
	where we used~\eqref{eq:cbc} for the first inequality and Lemma~\ref{lemma:stochastic_domination_nu} for the second one.
	
	Now, by Lemma~\ref{lemma:nu-convergence} and Proposition~\ref{prop:exp.decay_atrc01}, there exist~$\alpha, M > 0$ such that
	\[
		\nu_{\Omega}[0\leftrightarrow\partial\Lambda_n] \leq \mu_{J,U}^{0,1}[0\leftrightarrow\partial\Lambda_n]+e^{-\alpha n}
	\leq 8nMe^{-\alpha n}+e^{-\alpha n}.\qedhere
	\]
\end{proof}

It remains to show Lemmata~\ref{lemma:nu-convergence} and~\ref{lemma:nu-atrc}.

\begin{proof}[Proof of Lemma~\ref{lemma:nu-convergence}]
	By construction, $\nu_{\Omega}$ can be sampled from $\Spin_{\calD_{\Omega},c}^{+,+;c_b}$ using~\eqref{eq:ber-c-b}. 
	By Corollary~\ref{cor:omega-tau-from-spins}, $\mu_{\Omega,J,U}^{0,\mathbbm{1}}$ can be sampled from $\Spin_{\calD_{\Omega},c}^{+,+}$ using~\eqref{eq:ber-c}.
	The measures~$\Spin_{\calD_{\Omega},c}^{+,+;c_b}$ and~$\Spin_{\calD_{\Omega},c}^{+,+}$ both converge to~$\Spin_{c}^{+,+}$, as $\Omega\nearrow\bbL$, by Corollary~\ref{cor:spin-weak-limits}.
	Also, the rules~\eqref{eq:ber-c-b} and~\eqref{eq:ber-c} are local and coincide outside of the boundary (which is irrelevant in the limit).
	Thus, $\nu_{\Omega}$ and $\mu_{\Omega,J,U}^{0,\mathbbm{1}}$ have the same limit and it can be sampled from~$\Spin_{c}^{+,+}$ using the same rules.
	Their locality implies that the convergence inherits the exponential relaxation property~\eqref{eq:exp.relax._spin} (recall that~$\calD_{\Lambda_{4n}} = \Delta_{8n}$) and, by Proposition~\ref{prop:exp.decay_atrc01}, the limit is~$\mu_{J,U}^{0,1}$.
\end{proof}

\begin{proof}[Proof of Lemma~\ref{lemma:nu-atrc}]
	Fix a domain $\Omega$ in $\bbL$, and take $c_b:=e^{\lambda/2}$. 
	Recall that~$E_{\sigma^\bullet}$ denotes the set of disagreement edges in~$\sigma^\bullet$ (Section~\ref{sec:sixv-at}).
	
	\noindent\textbf{Step 1:} The measure~$\nu_\Omega$ can be written in the following form:
	\begin{equation}\label{eq:nu-formula}
	\nu_{\Omega}[\omega]\propto \left(\tfrac{2}{c-1}\right)^{\abs{\omega\setminus E_{\overline{\partial}\Omega}}}\left(\tfrac{2}{c_b-1}\right)^{\abs{\omega\cap E_{\overline{\partial}\Omega}}}2^{k(\omega)}\sum_{\substack{\sigma^\bullet\in\{\pm1\}^{\Omega}\\ \sigma^\bullet\perp\omega,\,\sigma^\bullet\vert_{\overline{\partial}\Omega}\equiv 1}}\left(\tfrac{1}{c-1}\right)^{\abs{E_{\sigma^\bullet}}}.
	\end{equation}
	For brevity, we write $E$ for $E_\Omega$ and $\partial E$ for $E_{\overline{\partial} \Omega}$. In a slight abuse of notation, we also set~$(E_{\sigma^\circ})^*=\{e^*:e\in E_{\sigma^\circ}\}$.
	The law of $(\sigma,\omega)$ defined by~\eqref{eq:ber-c-b} satisfies:
	\begin{align*}
		\prob[\sigma,\omega]&=&\, \Spin_{\calD_{\Omega},c}^{+,+;c_b}[\sigma]\,\ind{\sigma^\bullet\perp\omega,\,\sigma^\circ\perp\omega^*}
		&\times\left(\tfrac{1}{c}\right)^{\abs{\omega\setminus ((E_{\sigma^\circ})^* \cup \partial E)}}\left(\tfrac{c-1}{c}\right)^{\abs{E\setminus (\omega\cup E_{\sigma^\bullet} \cup \partial E)}}\\
		&& &\times\left(\tfrac{1}{c_b}\right)^{\abs{\omega\cap\partial E \setminus (E_{\sigma^\circ})^*}}
		\left(\tfrac{c_b-1}{c_b}\right)^{\abs{\partial E\setminus (\omega \cup E_{\sigma^\bullet})}}\\
		&\propto&\,(c-1)^{-\abs{(\omega \cup E_{\sigma^\bullet}) \setminus\partial E}}\,
		(c_b-&1)^{-\abs{(\omega \cup E_{\sigma^\bullet}) \cap\partial E}}
		\,\ind{\sigma\vert_{\partial\calD}\equiv 1}\,
		\ind{\sigma^\bullet\perp\omega,\,\sigma^\circ\perp\omega^*}.
	\end{align*}
	Note that $E_{\sigma^\bullet}\cap\partial E=\varnothing$ since $\sigma_{|\partial\Omega}\equiv 1$. 
	Summing over $\sigma$ then gives
	\begin{align*}
	\prob[\omega]\propto\,(c-1)^{-\abs{\omega\setminus\partial E}}(c_b-1)^{-\abs{\omega\cap\partial E}}\,2^{k^1(\omega^*)}\sum_{\substack{\sigma^\bullet\in\{\pm1\}^{\Omega}\\ \sigma^\bullet\perp\omega,\,\sigma^\bullet\vert_{\overline{\partial}\Omega}\equiv 1}}(c-1)^{-\abs{E_{\sigma^\bullet}}}.
	\end{align*}
	Finally, by Euler's formula (or induction), $k^1(\omega^*)=k(\omega)+\abs{\omega}+\mathrm{const}(\calD_{\Omega})$.
	
	\noindent\textbf{Step 2:} The marginal of $\atrc_{\Omega,{\rm w}_\tau,{\rm w}_{\tau\tau'}}^{0,\mathbbm{1}}$ on $\omega_\tau$ with weights~${\rm w}_\tau, {\rm w}_{\tau\tau'}$ given by~\eqref{eq:nu-atrc-weights} coincides with the right-hand side of~\eqref{eq:nu-formula}.
	
	Given $(\omega_{\tau},\omega_{\tau\tau'})\sim \atrc_{\Omega,{\rm w}_\tau,{\rm w}_{\tau\tau'}}^{0,\mathbbm{1}}$, define a spin configuration $\sigma^\bullet\in\{\pm 1\}^{V_\Omega}$ by assigning $+1$ to domain-boundary clusters of $\omega_{\tau\tau'}$ and $\pm1$ to interior clusters of~$\omega_{\tau\tau'}$ uniformly independently. Then their joint law can be written as
	\[
	\prob[\omega_\tau,\omega_{\tau\tau'},\sigma^\bullet]\propto \prod_{e\in \omega_\tau}{\rm w}_\tau(e)\cdot {\rm w}_{\tau\tau'}^{\abs{\omega_{\tau\tau'}\setminus\omega_\tau}}\,2^{k(\omega_\tau)}\,\ind{\omega_\tau\subseteq\omega_{\tau\tau'}}\,\ind{\sigma^\bullet\perp\omega_{\tau\tau'}}\,\ind{\sigma^\bullet\vert_{\overline{\partial} \Omega}\equiv 1}.
	\]
	Now, $\sigma^\bullet\perp\omega_{\tau\tau'}$ precisely if $\sigma^\bullet\perp\omega_\tau$ and $(\omega_{\tau\tau'}\setminus\omega_\tau)\cap E_{\sigma^\bullet}=\varnothing$. Sum over $\omega:=\omega_{\tau\tau'}\setminus\omega_\tau$:
	\[
	\prob[\omega_\tau,\sigma^\bullet]\propto\prod_{e\in\omega_\tau}{\rm w}_\tau(e) \cdot 2^{k(\omega_\tau)}\,\ind{\sigma^\bullet\perp\omega_{\tau}}\,\ind{\sigma^\bullet\vert_{\overline{\partial} \Omega}\equiv 1}\sum_{\omega\subseteq E_\Omega\setminus(\omega_\tau\cup E_{\sigma^\bullet})}{\rm w}_{\tau\tau'}^{\abs{\omega}}.
	\]
	The last term equals~$({\rm w}_{\tau\tau'}+1)^{\abs{E_\Omega}-\abs{\omega_\tau}-\abs{E_{\sigma^\bullet}}}$.
	Finally, summing over $\sigma^\bullet$, we arrive at 
	\begin{align}
	\prob[\omega_\tau]\propto\prod_{e\in \omega_\tau}\tfrac{{\rm w}_\tau(e)}{{\rm w}_{\tau\tau'}+1}\cdot 2^{k(\omega_\tau)}\sum_{\substack{\sigma^\bullet\in\{\pm 1\}^{V_\Omega}\\ \sigma^\bullet\perp\omega_\tau,\,\sigma^\bullet\vert_{\overline{\partial} \Omega}\equiv 1}}\left(\tfrac{1}{{\rm w}_{\tau\tau'}+1}\right)^{\abs{E_{\sigma^\bullet}}}.\label{eq:atrc_marginals_spins}
	\end{align}
	Plugging in the weights~\eqref{eq:nu-atrc-weights} while using that $\sinh 2J=e^{-2U}$ and $c=\coth(2J)$ gives that~\eqref{eq:atrc_marginals_spins} agrees with~\eqref{eq:nu-formula}, which finishes the proof. 
\end{proof}

\section{No infinite cluster in the wired self-dual ATRC}\label{sec:no-inf-cluster-1-1}

\begin{proposition}\label{prop:no-perco-atrc11}
	Let $0<J<U$ satisfy $\sinh 2J=e^{-2U}$. 
	Then, $\atrc_{J,U}^{1,1}[0\xlra{\omega_\tau}\infty]=0.$
\end{proposition}

The proof of Proposition~\ref{prop:no-perco-atrc11} again relies on the coupling with the six-vertex model, Proposition \ref{prop:6v-atrc-coupling}. First of all, by the non-coextistence theorem~\cite{She05,DumRaoTas19}, it is sufficient to show that $\atrc_{J,U}^{1,1}$ admits an infinite $\omega_\tau^*$-cluster. 
If the latter is not the case, the infinite-volume limit of the marginals of $\Spin_{\calD,c}^{+,\pm}$ on $\{\pm1\}^{\bbL^*}$ can be shown to be tail-trivial.
Exploring clusters of~$1$ and~$-1$ (in $\bbT$-connectivity) and using the non-coexistence theorem,
we obtain that the limit of $\HF_{\calD,c}^{0,\pm1}$ is either $\HF_c^{0,1}$ or $\HF_c^{0,-1}$, thereby contradicting the invariance of $\HF_{\calD,c}^{0,\pm1}$ under $h\mapsto -h$.

In the following remark, we summarise some basic properties of the ATRC marginals $\mu_{\Omega,J,U}^{1,1}$ and $\mu_{\Omega,J,U}^{\mathbbm{1},\mathbbm{1}}$ (defined in Section~\ref{sec:sixv-at}) and their infinite-volume limit that we will use in Sections~\ref{sec:inf-vol} and~\ref{sec:proof-no-perco-1-1}.

\begin{remark}\label{rem:atrc_basics}
Recall that, for domains $\Omega_k\nearrow\bbL$, the measures $\atrc_{\Omega_k,J,U}^{1,1}$ form a decreasing sequence and converge to~$\atrc_{J,U}^{1,1}$.
In particular, the same holds for the marginals on~$\omega_\tau$: $\mu_{\Omega_k,J,U}^{1,1}$ converges to $\mu_{J,U}^{1,1}$. Clearly $\mu_{\Omega_k,J,U}^{\mathbbm{1},\mathbbm{1}}$ converges to $\mu_{J,U}^{1,1}$ as well. 
It is then standard (\cite[Chapter 4.3]{Gri06}) that $\mu_{J,U}^{1,1}$ is invariant under translations and tail-trivial (and hence ergodic). 
Moreover, $\atrc_{\Omega_k,J,U}^{1,1}$ (and thus their limit and its marginals) satisfies the finite-energy property. 
Therefore, the Burton--Keane argument~\cite{BurKea89} and the non-coexistence theorem~\cite{She05,DumRaoTas19} apply.
\end{remark}

\subsection{Semi-free measures in infinite volume}
\label{sec:inf-vol}

In this section, we will show weak convergence for some finite-volume spin and height function measures defined in Section \ref{sec:sixv}.

\begin{lemma}\label{lem:spin_pm_inf_volume}
	Let~$0<J<U$ be on the self-dual line~\eqref{eq:self-dual} and take $c:=\coth 2J$. 
	Let~$\omega_\tau\sim \mu_{J,U}^{1,1}$.
	Define~$\chi_c^{+,\pm}$ as the distribution on~$\{\pm 1\}^{\bbL^*}$ obtained by assigning~$\pm 1$ to every cluster of~$\omega_\tau^*$ uniformly and independently.
	Then, for any sequence of odd domains $\calD_k\nearrow\bbZ^{2}$, the marginals of $\Spin_{\calD_k,c}^{+,\pm}$ on $\sigma^\circ$ converge weakly to~$\chi_c^{+,\pm}$.
	Moreover, $\chi_c^{+,\pm}$ is translation-invariant, positively associated and satisfies the finite-energy property.
\end{lemma}

\begin{proof}
	Fix~$J<U$. Let~$\calD_k$ be a sequence of odd domains on~$\bbZ^2$ and~$\Omega_k$ the corresponding subgraphs of~$\bbL$ such that~$\calD_k = \calD_{(\Omega_k)^*}$.
	Let~$\omega_\tau^k$ be sampled from~$\mu_{\Omega_k,J,U}^{\mathbbm{1},\mathbbm{1}}$.
	By Proposition \ref{prop:6v-atrc-coupling}, assigning~$\pm 1$ to clusters of $(\omega_\tau^k)^*$ uniformly independently gives the marginal of $\Spin_{\calD_k,c}^{+,\pm}$ on~$\sigma^\circ$.
	Since~$\mu_{\Omega_k,J,U}^{\mathbbm{1},\mathbbm{1}}$ converges to~$\mu_{J,U}^{1,1}$ that exhibits at most one infinite cluster in~$\omega_\tau^*$, the marginal of $\Spin_{\calD_k,c}^{+,\pm}$ on $\sigma^\circ$ converges to~$\chi_c^{+,\pm}$.
	
	Clearly, $\chi_c^{+,\pm}$ inherits translation-invariance and the finite-energy property from $\mu_{J,U}^{1,1}$. 
	By Proposition \ref{prop:spin-fkg}, the marginal of $\Spin_{\calD_k,c}^{+,\pm}$ on $\sigma^\circ$ satisfies the FKG inequality.
	Hence, the same holds for its limit~$\chi_c^{+,\pm}$.
\end{proof}

Working with measures on height functions (rather than spins) is more convenient as they satisfy stochastic ordering in boundary conditions. 
In the proof of Proposition~\ref{prop:no-perco-atrc11}, we use an infinite-volume version of $\HF_{\calD,c}^{0,\pm 1}$.
We show existence of such subsequential limit in the next lemma by sandwiching $\HF_{\calD,c}^{0,\pm 1}$ between $\HF_{\calD,c}^{0,-1}$ and~$\HF_{\calD,c}^{0,1}$.

\begin{lemma}\label{lem:hf-lim-pm-subseq}
	Let $c>2$. For any sequence of domains $\calD_k$ increasing to $\bbZ^{2}$, there exists a subsequence $(k_\ell)$ such that the measures $\HF_{\calD_{k_\ell},c}^{0,\pm1}$ converge weakly to some $\HF_c^{0,\pm 1}$ as $\ell$ tends to infinity.
\end{lemma}

\begin{remark}
	Proposition \ref{prop:perco_dual_atrc_11} and its proof allow to show that the limiting measure is~$\frac{1}{2}(\HF_c^{0,1}+\HF_c^{0,-1})$ for any (sub)sequence.
	We do not use this statement and omit the details.\end{remark}

\begin{proof}[Proof of Lemma \ref{lem:hf-lim-pm-subseq}]
	By \cite[Proposition 4.9]{Geo11}, it suffices to show that $(\HF_{\calD_k,c}^{0,\pm1})_{k\geq 1}$ is locally equicontinuous: for any finite $V \subset \bbZ^{2}$ and any decreasing sequence of local events $(A_m)_{m\geq 1}$ supported on $V$ and with $\cap_{m\geq 1}A_m=\varnothing$, it holds that 
	\[
	    \limsup_{k\to\infty}\HF_{\calD_k,c}^{0,\pm1}[A_m]\to 0 \quad\text{ as }m\to\infty. 
	\]
	By Proposition \ref{prop:hf-fkg}, finite-volume six-vertex height function measures are stochastically ordered with respect to the boundary conditions, whence
	\[
		\HF_{\calD_k,c}^{0,-1}\leq_{\mathrm{st}}\HF_{\calD_k,c}^{0,\pm1}\leq_{\mathrm{st}}\HF_{\calD_k,c}^{0,1}.
	\]
	Moreover, by Proposition~\ref{proposition:HF_weak_limit}, $\HF_{\calD_k,c}^{0,-1}$ converges to $\HF_c^{0,-1}$ and $\HF_{\calD_k,c}^{0,1}$ to $\HF_c^{0,1}$ as $k$ tends to infinity. These statements together easily imply the required local equicontinuity. 
\end{proof}

\subsection{Proof of Proposition \ref{prop:no-perco-atrc11}}
\label{sec:proof-no-perco-1-1}

As we argued in Remark~\ref{rem:atrc_basics}, it suffices to find a dual infinite cluster:

\begin{proposition}\label{prop:perco_dual_atrc_11}
	Let $0<J<U$ satisfy $\sinh 2J=e^{-2U}$. Then, $\atrc_{J,U}^{1,1}[0\xleftrightarrow{\omega_{\tau}^*}\infty]>0$.
\end{proposition}

\begin{remark}
	By the duality relation described in Section~\ref{sec:rc-at}, this is equivalent to saying that
	$\atrc_{J,U}^{0,0}[0\xleftrightarrow{\omega_{\tau\tau'}}\infty] >0$, for any $J<U$ on the self-dual line~\eqref{eq:self-dual}.
\end{remark}

The proof of Proposition \ref{prop:perco_dual_atrc_11} also relies on the non-coexistence theorem~-- but in the context of site percolation.
Following the notation of \cite{GlaPel19}, we let $\T^\circ$ be the graph with vertex set $\bbL^{*}$ where a vertex $(x,y)\in\bbL^{*}$ is adjacent to
\[(x,y)\pm (1,1),\,(x,y)\pm (1,-1)\text{ and }(x\pm 2,y).\]
Note that $\T^\circ$ is isomorphic to the triangular lattice. 

\begin{proof}[Proof of Proposition \ref{prop:perco_dual_atrc_11}]
	Fix~$J<U$. Recall that $\mu_{J,U}^{1,1}$ is the marginal of~$\atrc_{J,U}^{1,1}$ on~$\omega_\tau$.
	Assume for contradiction that $\mu_{J,U}^{1,1}$ does not admit an infinite dual cluster.
	
	Set $c=\coth 2J$. 
	Recall that~$\chi_c^{+,\pm}$ is obtained from $\omega\sim \mu_{J,U}^{1,1}$ by assigning uniformly independently~$\pm1$ to its dual clusters. 
	Since all them are finite by our assumption, $\chi_c^{+,\pm}$ inherits ergodicity from $\mu_{J,U}^{1,1}$.
	Also, by Lemma~\ref{lem:spin_pm_inf_volume}, $\chi_c^{+,\pm}$ is translation-invariant and satisfies the FKG inequality.
	Thus, by non-coexistence theorem, in~$\bbT^\circ$-connectivity, either $\chi_c^{+,\pm}$ admits no infinite cluster of minuses, whence
	\begin{equation} \label{eq:spin-pm-t-circuits}
		\chi_c^{+,\pm}(\exists \text{ infinitely many disjoint $\T^\circ$-circuits of $+$ around the origin})=1,
	\end{equation}
	or the same holds for $\T^\circ$-circuits of $-$.
	By symmetry, we can assume~\eqref{eq:spin-pm-t-circuits}.
		
	By Lemma \ref{lem:hf-lim-pm-subseq}, there exists a sequence of odd domains~$\calD_k$ such that $\HF_{\calD_k,c}^{0,\pm1}$ converge to some infinite-volume height function measure $\HF_c^{0,\pm 1}$ weakly. 
	Recall that~$\Spin_{\calD,c}^{+,\pm}$ is the push-forward of $\HF_{\calD,c}^{0,\pm1}$ under the modulo 4 mapping and, by Lemma~\ref{lem:spin_pm_inf_volume}, the marginals of~$\Spin_{\calD,c}^{+,\pm}$ on~$\sigma^\circ$ converge to~$\chi_c^{+,\pm}$ weakly.
	Then, \eqref{eq:spin-pm-t-circuits} implies that $\HF_c^{0,\pm1}$ admits infinitely many disjoint $\T^\circ$-circuits around the origin of constant height that is congruent to~$1$ modulo~$4$. 
	By Proposition \ref{prop:hf-fkg}, the measure $\HF_c^{0,\pm1}$ is between $\HF_c^{0,-1}$ and $\HF_c^{0,1}$ in the sense of stochastic domination.
	By Proposition~\ref{proposition:HF_weak_limit}, $\HF_c^{0,-1}$ and $\HF_c^{0,1}$ admit infinite clusters of $-1$ and $+1$, respectively.
	Hence, the above implies
	\begin{equation}\label{eq:hf-pm-circuits}
		\HF_c^{0,\pm1}(\exists \text{ infinitely many disjoint $\T^\circ$-circuits of $+1$ around the origin}) = 1.
	\end{equation}
	By a standard exploration argument and FKG inequality (details below), \eqref{eq:hf-pm-circuits} implies that~$\HF_c^{0,\pm1}$ stochastically dominates and hence equals $\HF_c^{0,1}$.
	This leads to a contradiction since~$\HF_c^{0,\pm1}$ is invariant under $h\mapsto -h$ while $\HF_c^{0,1}$ is not.
	
	It remains to show that~\eqref{eq:hf-pm-circuits} implies~$\HF_c^{0,\pm1}=\HF_c^{0,1}$.
	It is sufficient to prove~$\HF_c^{0,\pm1}[A]=\HF_c^{0,1}[A]$, for any increasing local event~$A$.
	Take any~$\eps>0$. 
	Since~$\HF_{\calD,c}^{0,1}$ converges to~$\HF_c^{0,1}$ weakly as~$\calD\nearrow\bbZ^2$,
	we can find~$n\geq 1$ such that, for any domain~$\calD$ containing~$\Delta_n$,
	\begin{equation}\label{eq:hf-limit-close}
		\vert\HF_{\calD,c}^{0,1}[A]-\HF_c^{0,1}[A]\vert<\eps.
	\end{equation}
	We can find $\calD\supseteq\Delta_n$ large enough such that
	\begin{equation}\label{eq:hf-t-circ-in-d}
		\left| \HF_c^{0,\pm1}[A]-\HF_c^{0,\pm1}[A\,|\, \exists \text{ $\bbT^\circ$-circuit in $\calD$ surrounding $\Delta_n$}]\right|<\eps.
	\end{equation}
	Let $\Cc$ be the outermost $\T^\circ$-circuit of height $+1$ in $\calD$ surrounding $\Delta_n$ (if such circuit does not exist, we set~$\Cc:=\emptyset$).
	By the domain Markov property and stochastic ordering in boundary conditions, for any~$\bbT^\circ$-circuit~$C$ surrounding~$\Delta_n$ and contained in~$\calD$,
	\begin{equation}\label{eq:hf-stoch-dom-circ}
		\HF_c^{0,\pm1}[A\,|\,\Cc=C]\geq\HF_{\calD_C,c}^{0,1}[A],
	\end{equation}
	where $\calD_C$ is the connected component of the origin in the graph obtained from~$\bbZ^2$ after removing all vertices on~$C$ or adjacent to it in~$\bbZ^2$.
	Since $\calD_C\supset\Delta_n$, by~\eqref{eq:hf-limit-close}, the right-hand side in~\eqref{eq:hf-stoch-dom-circ} is $\eps$-close to $\HF_c^{0,1}[A]$.
	Putting this together with~\eqref{eq:hf-t-circ-in-d}, we get $\HF_c^{0,\pm1}[A]\geq\HF_c^{0,1}[A]-2\eps$.
	Since~$\eps>0$ was arbitrary, we obtain~$\HF_c^{0,\pm1}[A]\geq\HF_c^{0,1}[A]$.
	The opposite inequality follows by the comparison of boundary conditions.
\end{proof}

\section{Proof of Proposition~\ref{thm:finite_expo_decay_at_sd}}
\label{sec:exp-dec-finite}

Our goal is to use exponential decay under~$0,1$ conditions in Proposition~\ref{prop:exp-decay-0-1} to improve non-percolation statement of Proposition~\ref{prop:no-perco-atrc11} and get exponential decay in finite volume stated in Proposition~\ref{thm:finite_expo_decay_at_sd}.
We use the approach of~\cite{CerMes10} and~\cite[Appendix]{CamIofVel08}.
Additional difficulties in our case come from a weaker domain Markov property of the ATRC measure.

Fix~$J<U$ and~$n\geq 1$.
For any vertex~$x\in \Lambda_n$, define
\begin{equation*}
\eta(x)=\mathbbm{1}_{\lbrace x \xlra{\omega_\tau}\partial\Lambda_{n}\rbrace}.
\end{equation*}

The next lemma provides a lower bound on the size of the boundary cluster of $\omega_{\tau}$.
Recall that we denote by~$\mu$ the marginal of the ATRC measure on~$\omega_\tau$.

\begin{lemma}\label{lemma:boundary_cluster_not_big}
For any $\delta>0$, there exists $\alpha:=\alpha(\delta,\betasd)>0$ such that 
\begin{equation*}
\mu^{1,1}_{\Lambda_{n},\beta_{sd}}
\left(\sum_{x\in\Lambda_{n}}\eta_{x}\geq\delta n^{2}\right)\leq e^{-\alpha n^{2}}.
\end{equation*}
\end{lemma}

\begin{proof}
	Fix $\delta>0$. 
	It follows from Proposition~\ref{prop:no-perco-atrc11} that 
	\begin{equation*}
	\lim_{n\rightarrow\infty}\mu^{1,1}_{\Lambda_{n},\betasd}(0\xlra{\omega_\tau}\partial \Lambda_{n})=0.
	\end{equation*}
	This implies that one can find $M:=M(\delta)>0$ such that
	\begin{equation*}
		\bbE_{\Lambda_{M},\betasd}^{1,1}\left[\tfrac{1}{\abs{\Lambda_{M}}}\sum_{x\in\Lambda_{M}}\eta_{x}\right]<\frac{\delta}{2}.
	\end{equation*}
	Fix $n\gg M$. Without loss of generality, assume that $n=(2k+1)M$. One has
	\begin{equation*}
	\dfrac{1}{\abs{\Lambda_{n}}}\sum_{x\in\Lambda_{n}}\eta_{x}
	\leq
	 \dfrac{1}{\abs{\Lambda_{k}}}\sum_{x\in\Lambda_{k}}
	 \left(\dfrac{1}{\abs{\Lambda_{M}}}\sum_{y\in 2Mx+\Lambda_{M}}
	 \mathbbm{1}_{y\leftrightarrow\partial (2Mx+\Lambda_{M})}\right).
	\end{equation*}
	Denote the expression in the brackets by~$Y_{x,M}$.
	Then,
	\begin{align*}
	\mu^{1,1}_{\Lambda_{n}, \betasd}\left(\dfrac{1}{\abs{\Lambda_{n}}}\sum_{x\in\Lambda_{n}}\eta_{x}\geq \delta \right)
	&\leq 
	\mu^{1,1}_{\Lambda_{n}, \betasd}
	\left(\dfrac{1}{\abs{\Lambda_{k}}}\sum_{x\in\Lambda_{k}}Y_{x, M}\geq\delta \right)
	 \\
	 &\leq 
	 \mu^{1,1}_{\Lambda_{n}, \betasd}\left(\dfrac{1}{\abs{\Lambda_{k}}}\sum_{x\in\Lambda_{k}}Y_{x, M}\geq \delta \, \Big| \, B_{M}\right),
	\end{align*}
	where $B_{M}=\bigcap_{x\in \Lambda_k}\{ \omega_{\tau}\vert_{\partial (2Mx+\Lambda_{M})} \equiv 1 \}$ and the last inequality uses~\eqref{eq:fkg} and that $B_{M}$ is increasing. 
	Note that under $\mu^{1,1}_{\Lambda_{n}, \betasd}(\cdot\vert\hphantom{,} B_{M})$, the random variables $Y_{x,M}$ are i.i.d. The statement then follows from Hoeffding's inequality.
\end{proof}

\begin{proof}[Proof of Proposition~\ref{thm:finite_expo_decay_at_sd}]
	We aim to show that, up to an arbitrary small exponential error, there exists a blocking surface of closed edges around $\Lambda_{\frac{4n}{5}}$ in $\Lambda_{n}$.
	
	For each~$\ell \in [1,n]$ and~$x\in \partial\Lambda_\ell$, define
	\[
		f(x):=\mathbbm{1}_{\{x\xlra{\omega_\tau, \Lambda_n\setminus \Lambda_\ell}\partial\Lambda_{n} \}}
		\hspace{5mm} \text{and} \hspace{5mm}
		N(\ell):= \sum_{x\in \partial\Lambda_\ell} f(x).
	\]
	Define~$A_\ell$ as the event that~$N(\ell) \leq \delta n$.
	Since~$f(x) \leq \eta(x)$, Lemma~\ref{lemma:boundary_cluster_not_big} implies that, up to an error~$e^{-cn^2}$, event~$A_\ell$ occurs for some~$\ell \in [4n/5,n]$, whence
	\begin{align*}
		\mu^{1,1}_{\Lambda_{n},\betasd}&(0\xlra{\omega_\tau}\partial\Lambda_{n/5}) 
		\leq \sfe^{-\alpha n^2} + \sum_{4n/5 \leq \ell \leq n} \mu^{1,1}_{\Lambda_{n}, \betasd}(0\xlra{\omega_\tau}\partial\Lambda_{n/5} \, \vert\, A_\ell) \cdot \mu^{1,1}_{\Lambda_{n},\betasd}(A_\ell)\\
		&\leq \sfe^{-\alpha n^2} + \sfe^{\alpha'\delta n}\sum_{4n/5 \leq \ell \leq n} 
		\mu^{1,1}_{\Lambda_{n},\betasd}(0\xlra{\omega_\tau}\partial\Lambda_{n/5}, \Lambda_{4n/5}\nxlra{\omega_\tau}\partial\Lambda_n \, \vert\, A_\ell) \cdot \mu^{1,1}_{\Lambda_{n}, \betasd}(A_\ell)\\
		&= \sfe^{-\alpha n^2} + \sfe^{\alpha'\delta n} \mu^{1,1}_{\Lambda_{n}, \betasd}(0\xlra{\omega_\tau}\partial\Lambda_{n/5}, \Lambda_{4n/5}\nxlra{\omega_\tau}\partial\Lambda_n),
	\end{align*}		
	where we used that~$A_\ell$ is measurable with respect to edges of~$\omega_\tau$ in~$\Lambda_n\setminus \Lambda_\ell$ and that, conditioned on~$A_\ell$, there are maximum~$4\delta n$ edges that are incident to vertices on~$\partial\Lambda_\ell$ that are connected to~$\Lambda_n$ and we can disconnect~$\Lambda_{4n/5}$ from~$\partial\Lambda_n$ by closing all these edges.

	On the event $\lbrace\Lambda_{4n/5}\nxlra{\omega_\tau}\partial\Lambda_n\rbrace$, there exists a circuit of closed edges in $\omega_{\tau}$ that surrounds~$\Lambda_{4n/5}$.
	Denote the exterior-most such circuit by $\zeta$ an explore it from the outside: 
	\begin{align*}
		\mu^{1,1}_{\Lambda_{n}, \betasd}(0\xlra{\omega_\tau} \partial\Lambda_{n/5}, \partial\Lambda_{4n/5}\nxlra{\omega_\tau}\partial\Lambda_{n})
		&= \sum_C\mu^{1,1}_{\Lambda_{n}, \betasd}(0\xlra{\omega_\tau}\partial\Lambda_{n/5} \, \vert \, \zeta=C) 
		\mu^{1,1}_{\Lambda_{n}, \betasd}( \zeta=C)\\
		&\leq \sum_C \mu^{0,1}_{\Omega_C, \betasd}(0\xlra{\omega_\tau}\partial\Lambda_{n/5})
		\mu^{1,1}_{\Lambda_{n}, \betasd}( \zeta=C),
	\end{align*}
	where the sum is over all possible values of~$\zeta$ and we define~$\Omega_C$ as the subgraph of~$\bbL$ bounded by~$C$; the inequality relies on~\eqref{eq:cbc} and on the~$0,1$ boundary conditions being domain Markov for~$\mu$.
	Note that~$\Omega_C$ can be turned into a domain by consecutively removing vertices of degree~$1$~-- denote it by~$\Omega_C'$.
	Such operations can only increase the measure, whence $\mu^{0,1}_{\Omega_C, \betasd} \leq_{{\rm st}} \mu^{0,1}_{\Omega_C', \betasd}$, and $\Lambda_{4n/5}\subset \Omega_C'$.
	Thus, the right-hand side in the last equation is exponentially small by Proposition~\ref{prop:exp-decay-0-1}.
	Combining the bounds, we get
	\begin{align*}
		\mu^{1,1}_{\Lambda_{n}, \betasd}(0\xlra{\omega_\tau}\partial\Lambda_{n/5})
		\leq \sfe^{-\alpha n^2} + \sfe^{\alpha'\delta n} \cdot \sfe^{-\alpha''n}.
	\end{align*}
	Taking~$\delta$ small enough finishes the proof.
\end{proof}

\section{The case $J\geq U$: Proof of Theorem~\ref{thm:sd-is-critical}}
\label{sec:0-le-u-le-j}

\subsection{ATRC for $J\geq U$}
\label{sec:atrc-0uj}

We fix $J\geq U$ and a finite subgraph $\Omega$ of $\bbL$.
The ATRC model is defined via an Edwards--Sokal-type expansion.
Since~$J\geq U$, the leading terms will correspond to interactions in~$\tau$ and in~$\tau'$.
Thus, the ATRC measure on~$\Omega$ with boundary conditions $\eta_\tau,\eta_{\tau'}$ is supported on pairs of percolation configurations $(\omega_\tau,\omega_{\tau'})\in\{0,1\}^{E_\Omega}\times\{0,1\}^{E_\Omega}$, and is defined by 
\begin{equation}\label{eq:def-atrc-0uj}
	\atrc^{\eta_\tau,\eta_{\tau'}}_{\Omega,J,U}(\omega_\tau, \omega_{\tau'})
	=\tfrac{1}{Z}\cdot 
	2^{k^{\eta_{\tau}}(\omega_{\tau})+k^{\eta_{\tau'}}(\omega_{\tau'})}\prod_{e\in E}a(\omega_\tau(e),\omega_{\tau'}(e)),
\end{equation}
where~$Z=Z(\Omega,J,U,\eta_\tau,\eta_{\tau'})$ is a normalizing constant and
\begin{equation}\label{eq:atrc-weights-0uj}
a(0,0):=e^{-4J}, \,
a(1,0) = a(0,1) :=e^{-2(J+U)}-e^{-4J}, \,
a(1,1):= 1-2e^{-2(J+U)} + e^{-4J}.
\end{equation}
Similarly to~\eqref{eq:atrc-repr.}, if $J>U$, we can write the measure as
\begin{equation}\label{eq:atrc-repr-0uj}
	\atrc^{\eta_\tau,\eta_{\tau'}}_{\Omega,J,U}(\omega_\tau, \omega_{\tau'})\propto \left(\tfrac{a(1,0)}{a(0,0)}\right)^{\abs{\omega_\tau}+\abs{\omega_{\tau'}}}\left(\tfrac{a(0,0)a(1,1)}{a(1,0)^2}\right)^{\abs{\omega_\tau\cap\omega_{\tau'}}}2^{k^{\eta_{\tau}}(\omega_{\tau})+k^{\eta_{\tau'}}(\omega_{\tau'})}.
\end{equation}
As before, if the parameters $J,U$ are fixed, we write $\atrc^{\eta_\tau,\eta_{\tau'}}_{\Omega,\beta}$ for $\atrc^{\eta_\tau,\eta_{\tau'}}_{\Omega,\beta J,\beta U}$.

\paragraph{Coupling of ATRC and AT.}
As mentioned above, edges in~$\omega_\tau$ and in~$\omega_{\tau'}$ describe interactions in~$\tau$ and in~$\tau'$.
In contrast to~\eqref{eq:coupling}, the correlations of the product $\tau\tau'$ are described by  simultaneous connections in both~$\omega_\tau$ and~$\omega_{\tau'}$: for any vertex $x\in V_\Omega$,
\begin{equation}\label{eq:atrc-at-coupling-0uj}
\langle\tau_x\rangle_{\Omega,J,U}^{+,+}=\atrc_{\Omega,J,U}^{1,1}(x\xlra{\omega_\tau}\partial\Omega),
\quad
\langle\tau_x\tau'_x\rangle_{\Omega,J,U}^{+,+}=\atrc_{\Omega,J,U}^{1,1}(x\xlra{\omega_\tau}\partial\Omega,x\xlra{\omega_{\tau'}}\partial\Omega).
\end{equation}
The statement extends to infinite volume in a standard way, see~\cite[Proposition 5.11]{Gri06}.

\paragraph{Basic properties.} 
The analogues of the properties~\eqref{eq:fkg}, \eqref{eq:cbc}, \eqref{eq:mon}, \eqref{eq:strong-mon} and~\eqref{eq:dmp} hold in this context as well. In particular, the measures $\atrc_{\Omega,J,U}^{0,0}$ and $\atrc_{\Omega,J,U}^{1,1}$ converge weakly to some $\atrc_{J,U}^{0,0}$ and $\atrc_{J,U}^{1,1}$, respectively, as $\Omega\nearrow\bbL$. 

\paragraph{Duality.} 
Given an ATRC configuration $(\omega_\tau,\omega_{\tau'})$, we define its dual $(\hat{\omega}_\tau,\hat{\omega}_{\tau'}) := (\omega_\tau^*,\omega_{\tau'}^*)$.
This extends Lemma~\ref{lemma:duality_relation} to all~$J,U>0$.

\subsection{Proof of Theorem~\ref{thm:sd-is-critical}}
Fix $J\geq U$. 
By~\eqref{eq:easy-betac-ineq}, $\betac^{\tau}\geq \betac^{\tau\tau'}$.
The opposite inequality follows directly from the ATRC representation and the coupling~\eqref{eq:atrc-at-coupling-0uj}. Indeed, $\betac^{\tau\tau'}\geq\betac^{\tau}$ since
\begin{equation*}
	\langle \tau_0\tau_0' \rangle_{\beta}^{+,+} 
	= \atrc_\beta^{1,1}(0\xlra{\omega_\tau} \infty, 0\xlra{\omega_{\tau'}} \infty)
	\leq \atrc_\beta^{1,1}(0\xlra{\omega_\tau} \infty)
	= \langle \tau_0 \rangle_{\beta}^{+,+}.
\end{equation*}
So we get~$\betac^{\tau}=\betac^{\tau\tau'}=:\betac$.
Similarly, $\betac^{\tau,f}=\betac^{\tau\tau',f}=:\betac^f$.

To derive part (ii) of Theorem~\ref{thm:sd-is-critical}, it remains to show that $\betac=\betasd$. The inequality $\betac\leq\betasd$ follows from a standard argument once the transition is shown to be \emph{sharp}: for $\beta<\betac$, there exists $c=c(\beta)>0$ such that, for every $n\geq1$,
\begin{equation}\label{eq:sharpness-0uj}
\atrc_\beta^{1,1}(0\xlra{\omega_\tau}\partial\Lambda_n)\leq e^{-cn}.
\end{equation}
This can be derived via a general approach~\cite{DumRaoTas19}, see Appendix~\ref{sec:sharpness}.

The reverse inequality is a consequence of Zhang's argument provided that~$\betac=\betac^f$, i.e. the transitions for the free and wired measures occur at the same point.
This follows from an analogue of Lemma~\ref{lemma:dense_subsets} (see Appendix~\ref{sec:dense_subsets} for the proof of both lemmata):

\begin{lemma}\label{lemma:f-equals-w-0uj}
	There exists~$D\subseteq\{(J,U)\in\bbR^2:J\geq U > 0\}$ with Lebesgue measure~$0$ such that, for any~$(J,U)\in D^c$, one has
	\begin{equation*}
		\atrc_{J,U}^{0,0}=\atrc_{J,U}^{1,1}.
	\end{equation*} 
\end{lemma}

\begin{proof}[Proof of Theorem~\ref{thm:sd-is-critical}]
	Fix $J\geq U$. Part (i) follows from Lemma~\ref{lemma:f-equals-w-0uj} and~\eqref{eq:strong-mon} in the same way as for $J<U$, see Section~\ref{sec:duality-critical-lines}.
	
	Recall the definition of the event $\calH_n^\tau$ in Section~\ref{sec:duality-critical-lines}. 
	By duality, symmetry and~\eqref{eq:cbc},
	\begin{equation}\label{eq:calh-sum}
		\atrc_{\betasd}^{0,0}(\calH_n^\tau)\leq\frac{1}{2}\leq\atrc_{\betasd}^{1,1}(\calH_n^\tau).
	\end{equation}
	If~$\betac>\betasd$, then, by~\eqref{eq:sharpness-0uj}, $\atrc_{\betasd}^{1,1}(\calH_n^\tau)$ converges to~$0$ as $n$ tends to infinity, 
	If~$\betac<\betasd$, then (since $\betac=\betac^f$) Zhang's argument implies that $\atrc_{\betasd}^{0,0}(\calH_n^\tau)$ converges to 1 as $n$ tends to infinity.
	Both statements contradict~\eqref{eq:calh-sum}.
\end{proof}

\appendix

\section{Sharpness}
\label{sec:sharpness}

The proof of sharpness for FK-percolation via the OSSS inequality~\cite{DumRaoTas19,OdoSakSchSer05} adapts to the ATRC.
For completeness, we present a sketch of this argument and give details for the steps that are specific for the ATRC.

\subsection{Sharpness for $J<U$}

We start by bounding the derivative in~$\beta$ by a covariance:

\begin{lemma}\label{lemma:eq_derivative}
	Let~$0<J<U$ and $\varepsilon>0$.
	Then, there exists $c=c(\varepsilon,J,U)>0$ such that, for any finite subgraph~$\Omega\subseteq\bbL$, any increasing event $A$, and any~$\beta_0\in [\varepsilon,\varepsilon^{-1}]$,
	\begin{equation*}
		\left(\frac{d}{d\beta}\atrc_{\Omega,\beta}^{1,1}[A]\right)\Bigg\vert_{\beta=\beta_0}\geq c\sum\limits_{e\in E_{\Omega}}\Cov[\1_{A},\omega_{\tau}(e)]+\Cov[\1_{A},\omega_{\tau\tau'}(e)].
	\end{equation*}
\end{lemma}

\begin{proof}
	Fix $J<U$. Recall that we can write the measure $\atrc^{1,1}_{\Omega,\beta}$ as in~\eqref{eq:atrc-repr.} with weights given by~\eqref{eq:atrc-weights-2}.
	Then, as in FK-percolation, we get a covariance formula:
	\begin{equation*}
		\dfrac{d}{d\beta}\atrc^{1,1}_{\Omega,\beta}[A]
		=\Cov\left[\1_{A},
		\dfrac{{\rm w}_\tau'(\beta)}{{\rm w}_\tau(\beta)}\cdot \abs{\omega_{\tau}}+
		\dfrac{{\rm w}_{\tau\tau'}'(\beta)}{{\rm w}_{\tau\tau'}(\beta)}\cdot \abs{\omega_{\tau\tau'}\setminus\omega_\tau}\right].
	\end{equation*}
	Since~$J<U$, we have ${\rm w}_{\tau\tau'}'(\beta), {\rm w}_{\tau}'(\beta) > 0$ and the statement follows.
\end{proof} 

Fix~$J<U$.
We prove sharpness only for~$\omega_\tau$, since the proof for~$\omega_{\tau\tau'}$ is the same. 
Recall that $\mu_{\Omega,\beta}^{1,1}$ is the marginal of $\atrc^{1,1}_{\Omega,\beta}$ on~$\omega_\tau$.
The key step in the proof of sharpness in~\cite{DumRaoTas19} is the extension of the OSSS inequality~\cite{OdoSakSchSer05} to dependent measures.
The inequality holds for any monotonic (FKG) measure on~$\{0,1\}^E$, for a finite set of edges~$E$.
In particular, it applies also to~$\mu^{1,1}_{\Lambda_{2n}, \beta}$, for~$n\geq 1$.
Instead of stating the OSSS inequality, we state its consequence that can be derived in the same way as in~\cite{DumRaoTas19}:

\begin{lemma}[\cite{DumRaoTas19}, Lemma 3.2]\label{lem:sharpness-cov-bound}
	For any $n\geq 1$, one has 
	\begin{equation*}
		\sum\limits_{e\in E_{2n}}\Cov[\1_{\{0\leftrightarrow\partial\Lambda_{n}\}},\omega_{e}]
		\geq
		\dfrac{n}{16\sum\limits_{k=0}^{n-1}\mu^{1,1}_{\Lambda_{2k},\beta}[0\leftrightarrow\partial\Lambda_{k}]}\mu^{1,1}_{\Lambda_{2n},\beta}[0\leftrightarrow\partial\Lambda_{n}]
		\left(1-\mu^{1,1}_{\Lambda_{2n},\beta}[0\leftrightarrow\partial\Lambda_{n}]\right),
	\end{equation*}
	where the covariance is taken with respect to the measure $\mu^{1,1}_{\Lambda_{2n},\beta}$.
\end{lemma}

We proceed as in~\cite{DumRaoTas19}.
Fix $\beta_{0}>0$. For $n,k\geq 1, \varepsilon<1$ and $\beta\in [\varepsilon,\varepsilon^{-1}]$, define
\begin{equation*}
\theta_{k}(\beta):=\mu^{1,1}_{\Lambda_{2k},\beta}[0\leftrightarrow\partial\Lambda_{k}]
\hspace{1.5cm}
S_{n}:=\sum_{k=0}^{n-1}\theta_{k}.
\end{equation*}
Lemma~\ref{lemma:eq_derivative} applied to 
$A=\lbrace 0\overset{\omega_{\tau}}{\longleftrightarrow}\partial\Lambda_{n}\rbrace$ and $\Omega=\Lambda_{2n}$ implies
\begin{equation}\label{equation:FKG_necessary}
\theta^{'}_{n}(\beta)\geq c\sum\limits_{e\in E_{2n}}\Cov[\1_{A},\omega_{\tau\tau'}(e)]+\Cov[\1_{A},\omega_{\tau}(e)]
\geq 
c
\sum\limits_{e\in E_{2n}}\Cov[\1_{A},\omega_{\tau}(e)],
\end{equation}
where we used FKG inequality in the last line. By Lemma~\ref{lem:sharpness-cov-bound},
\begin{equation*}
\theta^{'}_{n}
\geq 
c\dfrac{n}{16S_{n}}\theta_{n}(1-\theta_{n}).
\end{equation*}
By~\eqref{eq:cbc} and monotonicty in~$\beta \leq \varepsilon^{-1}$, we have $\theta_{n}(\beta)\leq\theta_{1}(\varepsilon^{-1})$.
Then, for some~$c_1>0$,
\begin{equation}\label{eq:differentialinequality}
\theta^{'}_{n}
\geq c_{1}\dfrac{n}{S_{n}}\theta_{n}.
\end{equation}
By~\cite[Lemma 3.1]{DumRaoTas19}, this inequality implies sharpness of the phase transition.

\subsection{Sharpness for~$J\geq U$}
The proof is the same as for $J < U$ and we only show the analogue of Lemma~\ref{lemma:eq_derivative}.

\begin{lemma}\label{lemma:eq-derivative-0uj}
Let $J\geq U > 0$ and $\eps>0$. Then, there exists $c=c(\eps,J,U)>0$ such that, for any finite subgraph $\Omega\subseteq\bbL$, any increasing event $A$, and any $\beta_0\in[\eps,\eps^{-1}]$,
\begin{equation*}
	\left(\frac{d}{d\beta}\atrc_{\Omega,\beta}^{1,1}[A]\right)\Bigg\vert_{\beta=\beta_0}\geq c\sum_{e\in E_\Omega} \left(\Cov[\1_A,\omega_\tau(e)]+\Cov[\1_A,\omega_{\tau'}(e)]\right). 
\end{equation*}
\end{lemma}
\begin{proof}
	For $J=U$, the model reduces to FK-percolation with cluster-weight~$q=4$, and the statement follows from~\cite[Theorem 3.12]{Gri06}.
	Fix $J>U$. Define $r(\beta),s(\beta)>0$ by
	\[
		r(\beta)= \tfrac{a(1,0)}{a(0,0)} = \tfrac{a(0,1)}{a(0,0)}  \quad\text{and}\quad s(\beta)=\tfrac{a(0,0)a(1,1)}{a(1,0)^2},
	\]
	where the $a(i,j)$ are given by~\eqref{eq:atrc-weights-0uj} evaluated at $(\beta J,\beta U)$.
	Recall that we can write the ATRC measure as in~\eqref{eq:atrc-repr-0uj}.
	Then, for any $c>0$ and any increasing event $A$, 
	\[
		\frac{d}{d\beta}\atrc_{\Omega,\beta}^{1,1}[A]-c\sum_{e\in E_\Omega}\Cov[\1_A,\omega_\tau(e)+\omega_{\tau'}(e)]=\Cov[\1_A,X_c],
	\]
	where 
	\[
		X_c=\sum_{e\in E_\Omega}\left(\tfrac{r'}{r}-c\right)(\omega_\tau(e)+\omega_{\tau'}(e))+\tfrac{s'}{s}\,\omega_\tau(e)\,\omega_{\tau'}(e).	
	\]
	It is easy to see that~$X_c$ is increasing in~$\omega$ when~$\beta\in[\eps,\eps^{-1}]$ and $c$ is small enough.
	Then, by~\eqref{eq:fkg}, $\Cov[\1_A,X_c]\geq 0$ and this finishes the proof.
\end{proof}

\section{Proof of Lemma \ref{lemma:exp.relax.rc_pc}}
\label{sec:exp-rel-fk}

\begin{proof}[Proof of Lemma \ref{lemma:exp.relax.rc_pc}]
	Fix $q>4$ and $p=p_{\rm sd}(q)$ and omit them in the notation below.
	Let $n\geq 1$ and $\Omega \supset \Lambda_{2n}$ be a finite subgraph of $\bbL$.
	By Proposition~\ref{prop:rc_basic_properties} and Strassen's theorem, there exists a coupling $\prob$ of $\eta_-\sim\FK^{\mathrm{w}}$ and $\eta_+\sim\FK_{\Omega}$ such that $\prob(\eta_-\leq\eta_+) = 1$.
	
	Define $\Cc$ to be the outermost circuit of edges in~$\eta_-$ surrounding $\Lambda_n$ and contained in $\Lambda_{2n}$ (if there is no such circuit, set~$\Cc:=\emptyset$).
	By exponential decay of connections for the dual $\eta_-^*f$ (Theorem~\ref{Thm:critical_rc_properties}), there exists $\alpha>0$ such that, for any~$n\geq 1$,
	\[
		\Prob{\Cc=\varnothing}\leq \FK^{\mathrm{w}}[\Lambda_{n}^*\xlra{\eta_-^*} \Lambda_{2n}^*]<8ne^{-\alpha n}. 
	\]
	Take any event $A$ depending only on edges in $\Lambda_{n}$. We have
	\begin{align*}
	    \FK_\Omega^{\mathrm{w}}[A]-\FK^{\mathrm{w}}[A]
	    &\leq \Prob{\eta_+\in A,\eta_-\not\in A, \Cc\neq\varnothing}+8ne^{-\alpha n}\\
	    &=\sum_{C \neq \emptyset}\Prob{\eta_+\in A,\eta_-\not\in A\mid\Cc=C}\Prob{\Cc=C}+8ne^{-\alpha n} = 8ne^{-\alpha n},
	\end{align*}
	where the sum runs over all realisations $C\neq \emptyset$ of $\Cc$.
	Above we used that~$\{\Cc = C\}$ is measurable with respect to the exterior of~$C$; 
	on~$\{\Cc = C\}$, the distributions of~$\eta_+$ and~$\eta_-$ in the interior of~$C$ are equal;
	since~$\prob(\eta_-\leq \eta_+)=1$, the latter implies~$\1_A(\eta_-)=\1_A(\eta_+)$.
\end{proof}

\section{Proof of Lemma~\ref{lemma:stoch.dom.atrc}}
\label{sec:proof-stoch-dom-atrc}

\begin{proof}[Proof of Lemma~\ref{lemma:stoch.dom.atrc}]
	By a generalization of the Holley criterion, see \cite[Section 4]{GeoHagMae01}, it suffices to check that, for all $e\in E_\Omega$ and~$(\zeta_\tau,\zeta_{\tau\tau'})\in\{(0,0),(0,1),(1,1)\}^{E_\Omega\setminus\{e\}}$, both
	\begin{align*}   
		f_e({\rm w}_\tau,\zeta_\tau,\zeta_{\tau\tau'})&:=\atrc_{\Omega,{\rm w}_\tau,{\rm w}_{\tau\tau'}}^{0,\mathbbm{1}}\left[\omega_\tau(e)=1\,\Big|\,(\omega_\tau,\omega_{\tau\tau'})\vert_{E_\Omega\setminus\{e\}}=(\zeta_\tau,\zeta_{\tau\tau'})\right],\\
		g_e({\rm w}_\tau,\zeta_\tau,\zeta_{\tau\tau'})&:=\atrc_{\Omega,{\rm w}_\tau,{\rm w}_{\tau\tau'}}^{0,\mathbbm{1}}\left[\omega_{\tau\tau'}(e)=1\,\Big|\,(\omega_\tau,\omega_{\tau\tau'})\vert_{E_\Omega\setminus\{e\}}=(\zeta_\tau,\zeta_{\tau\tau'})\right]   
	\end{align*}
	are increasing in $(\zeta_\tau,\zeta_{\tau\tau'})$ and ${\rm w}_{\tau}$. 
	
	Fix $e,\zeta_\tau,\zeta_{\tau\tau'}$ as above. 
	Write $\underline{\zeta_{\tau}}$ and $\overline{\zeta_{\tau}}$ for the configurations that agree with $\zeta_\tau$ on $E_\Omega\setminus\{e\}$ while $\underline{\zeta_{\tau}}(e)=0$ and $\overline{\zeta_{\tau}}(e)=1$, and analogously for $\zeta_{\tau\tau'}$. Then,
	\begin{align*}
	    f_e({\rm w}_\tau,\psi,\zeta)
	    &=\frac{{\rm w}_\tau(e)2^{k(\overline{\zeta_\tau})+k^1(\overline{\zeta_{\tau\tau'}})}}{2^{k(\underline{\zeta_\tau})+k^1(\underline{\zeta_{\tau\tau'}})}+{\rm w}_{\tau\tau'}(e)2^{k(\underline{\zeta_\tau})+k^1(\overline{\zeta_{\tau\tau'}})}+{\rm w}_\tau(e)2^{k(\overline{\zeta_\tau})+k^1(\overline{\zeta_{\tau\tau'}})}}\\
	    &=\left(\frac{2^{k(\underline{\zeta_\tau})-k(\overline{\zeta_\tau})}(2^{k^1(\underline{\zeta_{\tau\tau'}})-k^1(\overline{\zeta_{\tau\tau'}})}+{\rm w}_{\tau\tau'}(e))}{{\rm w}_\tau(e)}+1\right)^{-1},
	\end{align*}
	which is clearly increasing in ${\rm w}_\tau$. Moreover $k(\underline{\zeta_\tau})-k(\overline{\zeta_\tau})$ and $k^1(\underline{\zeta_{\tau\tau'}})-k^1(\overline{\zeta_{\tau\tau'}})$ are decreasing in $(\zeta_\tau,\zeta_{\tau\tau'})$, respectively. 
	
	We now check that $1-g_e$ is decreasing in ${\rm w}_\tau$ and~$(\zeta_\tau,\zeta_{\tau\tau'})$. We have 
	\begin{align*}
	    1-g_e({\rm w}_\tau,\zeta_\tau,\zeta_{\tau\tau'})
	    &=\frac{2^{k(\underline{\zeta_\tau})+k^1(\underline{\zeta_{\tau\tau'}})}}{2^{k(\underline{\zeta_\tau})+k^1(\underline{\zeta_{\tau\tau'}})}+{\rm w}_{\tau\tau'}(e)2^{k(\underline{\zeta_\tau})+k^1(\overline{\zeta_{\tau\tau'}})}+{\rm w}_\tau(e)2^{k(\overline{\zeta_\tau})+k^1(\overline{\zeta_{\tau\tau'}})}}\\
	    &=\left(1+2^{k^1(\overline{\zeta_{\tau\tau'}})-k^1(\underline{\zeta_{\tau\tau'}})}({\rm w}_{\tau\tau'}(e)+{\rm w}_\tau(e)2^{k(\overline{\zeta_\tau})-k(\underline{\zeta_\tau})})\right)^{-1},
	\end{align*}
	which is clearly decreasing in ${\rm w}_\tau$ and~$(\zeta_\tau,\zeta_{\tau\tau'})$.
\end{proof}

\section{Equality of infinite volume $\atrc$ measures}
\label{sec:dense_subsets}
In this section, we prove Lemmata~\ref{lemma:dense_subsets} and~\ref{lemma:f-equals-w-0uj}. 

\paragraph{The case $J<U$.}
The following statement will imply Lemma~\ref{lemma:dense_subsets}.

\begin{lemma}\label{lemma:atrc_00vs11}
	There exist two families of smooth curves $(\gamma^\tau_r)_{r>0}$ and $(\gamma^{\tau\tau'}_s)_{s>0}$ such that~$(i)$: for any~$0<J<U$, there exist~$s,r>0$ such that~$(J,U)\in \gamma_r^\tau\cap \gamma_s^{\tau\tau'}$, and~$(ii)$: for any~$r,s>0$, the number of points~$(J,U)$ on~$\gamma_r^\tau$ (resp.~$\gamma_s^{\tau\tau'}$) such that the marginals of~$\atrc_{J,U}^{0,0}$ and~$\atrc_{J,U}^{1,1}$ on~$\omega_\tau$ (resp.~$\omega_{\tau\tau'}$) differ is at most countable.
\end{lemma}

This implies that the set of pairs $J<U$, for which~$\atrc_{J,U}^{0,0}$ and~$\atrc_{J,U}^{1,1}$ have different marginals on $\omega_\tau$ or $\omega_{\tau\tau'}$, has Lebesgue measure~$0$. 
By a monotone coupling argument, equality of both marginals implies 
that~$\atrc_{J,U}^{0,0} = \atrc_{J,U}^{1,1}$, and Lemma~\ref{lemma:dense_subsets} follows. 
We mention that $(\gamma^\tau_r)_{r>0}$ and~$(\gamma^{\tau\tau'}_s)_{s>0}$ are dual to each other.

\begin{proof}[Proof of Lemma \ref{lemma:atrc_00vs11}]
We follow a strategy presented in \cite[Theorem 1.12]{Dum17a} which is a rephrased version of an argument in \cite{LebMar72}. For any~$J<U$, any finite subgraph $\Omega\subseteq\bbL$ and any boundary conditions $\eta_\tau$ and $\eta_{\tau\tau'}$, we can write 
\[
	\atrc_{\Omega,J,U}^{\eta_\tau,\eta_{\tau\tau'}}[\omega_\tau,\omega_{\tau\tau'}]\propto \left(\tfrac{a(0,0)}{a(0,1)}\right)^{\abs{E_\Omega}-\abs{\omega_{\tau\tau'}}}\left(\tfrac{a(1,1)}{a(0,1)}\right)^{\abs{\omega_\tau}}2^{k^{\eta_{\tau}}(\omega_{\tau})+k^{\eta_{\tau\tau'}}(\omega_{\tau\tau'})}\,\ind{\omega_\tau\subseteq\omega_{\tau\tau'}}.
\]
Consider the curves where $a(0,0)/a(0,1)$ is constant (precisely if $U-J$ is constant). Fix an edge $e$ of $\bbL$. The Holley criterion~\cite{Hol74} easily gives that the function $(J,U)\mapsto\atrc_{J,U}^{1,1}[\omega_\tau(e)]$ is increasing along these curves, see e.g. the proof of~\cite[Lemma 11.14]{Gri06}. In particular, the set of discontinuity points is countable along each of them. Fix $C>0$ and $(J,U),(J',U')$ on the curve $a(0,0)/a(0,1)\equiv C$ with $J' < J$ (then also $U' < U$), and assume that $(J,U)$ is a continuity point. Define $a=\atrc_{J,U}^{0,0}[\omega_\tau(e)]$ and $b=\atrc_{J',U'}^{1,1}[\omega_\tau(e)]$. 

Note that, along the curve $a(0,0)/a(0,1)\equiv C$, the quantity $a(1,1)/a(0,1)$ is strictly increasing. Using this fact, analogous reasoning as in \cite{Dum17a} gives that $a\geq b$. Letting $(J',U')$ tend to $(J,U)$ along $a(0,0)/a(0,1)\equiv C$ and using that $(J,U)$ is a continuity point of $(J',U')\mapsto b$, we deduce 
\[
\atrc_{J,U}^{0,0}[\omega_\tau(e)]\geq \atrc_{J,U}^{1,1}[\omega_\tau(e)].
\]
Since the reversed inequality follows from \eqref{eq:cbc}, we obtain equality. By considering a monotone coupling of the corresponding marginals on $\omega_\tau$, it is easy to see that this implies that the marginals on $\omega_\tau$ coincide, and (i) is proved. The second statement (ii) is derived along the same lines when considering the curves where $a(1,1)/a(0,1)$ is constant. 

\end{proof}

\paragraph{The case $J\geq U$.}
Below is the analogue of Lemma~\ref{lemma:atrc_00vs11} that implies Lemma~\ref{lemma:f-equals-w-0uj}.
\begin{lemma}
	There exists a family of smooth curves $(\gamma_r^\tau)_{r>0}$ such that~$(i)$: for any~$J\geq U>0$, there exists~$r>0$ for which~$(J,U)\in \gamma_r^\tau$, and~$(ii)$: for any $r>0$, there exist only countably many points $(J,U)$ on $\gamma_r^\tau$ for which~$\atrc_{J,U}^{0,0}\neq\atrc_{J,U}^{1,1}$.
\end{lemma}

\begin{proof}
	For $J=U$, the model reduces to FK-percolation with cluster-weight $q=4$, and the statement follows from~\cite[Theorem 1.12]{Dum17a}. 
	Recall that, for any~$J > U$ and any finite subgraph $\Omega\subseteq\bbL$, we can write the measure as in~\eqref{eq:atrc-repr-0uj}. Consider the curves where the weight ${\rm w}:=(a(0,0)a(1,1))/(a(1,0)^2)$ is constant. The Holley criterion again allows to show that the ATRC measures are stochastically ordered along these lines. Moreover, the weight $a(1,0)/a(0,0)$ is increasing along each of them.
	
	Fix an edge $e$ of $\bbL$ and $C\geq 1$, and let $(J,U)$ be a continuity point of $(J',U')\mapsto \atrc_{J',U'}^{1,1}[\omega_\tau(e)]$ along the curve ${\rm w} \equiv C$. Take $(J',U')$ on the same curve with $J'<J$. Analogous reasoning as in the proof of~\cite[Theorem 1.12]{Dum17a} gives 
	\[
	\atrc_{J',U'}^{1,1}[\omega_\tau(e)]\leq\atrc_{J,U}^{0,0}[\omega_\tau(e)],
	\]
	with the only difference that one has to apply~\eqref{eq:fkg} to control $\abs{\omega_\tau}$ and $\abs{\omega_{\tau'}}$ simultaneously. Letting $(J',U')$ tend to $(J,U)$ along ${\rm w}\equiv C$ from below and using~\eqref{eq:cbc} gives $\atrc_{J,U}^{0,0}[\omega_\tau(e)]=\atrc_{J,U}^{1,1}[\omega_\tau(e)]$. A monotone coupling argument and symmetry between $\omega_\tau$ and $\omega_{\tau'}$ finish the proof.
\end{proof}

\bibliography{biblicomplete}
\bibliographystyle{plain}

\end{document}